\documentclass[11pt]{article}
\pdfoutput=1

\usepackage{lineno,hyperref}
\usepackage{color}
\usepackage{amsfonts}
\usepackage{amsmath}
\usepackage{amssymb}
\usepackage{amsthm}
\usepackage{subeqnarray}
\usepackage{lscape,graphicx,url}
\usepackage{anysize}
\usepackage{proof}
\usepackage{fancyhdr}
\usepackage{latexsym}
\usepackage{array}
\usepackage{multirow}
\usepackage{textcomp}
\usepackage{optprog}
\usepackage{mathtools}
\usepackage{proof}
\usepackage{natbib}
\usepackage[margin=1.5cm]{geometry}

\usepackage{enumitem}
\usepackage{csquotes}
\usepackage{comment}
\usepackage{bm}
\usepackage{checkend}
\usepackage{array}
\newcolumntype{H}{>{\setbox0=\hbox\bgroup}c<{\egroup}@{}}

\newcommand{\GG}[1]{}

\title{Order batching using an approximation for  the distance travelled by pickers}
\author{Cristiano Arbex Valle$^1$ \and John E Beasley$^2$}
\date{}

\begin{document}

\maketitle
\begin{center} 
{\footnotesize

$^1$Departamento de Ci\^{e}ncia da Computa\c{c}\~{a}o, \\
 Universidade Federal de Minas Gerais, \\
 Belo Horizonte, MG 31270-010, Brasil \\
 arbex@dcc.ufmg.br \\ \vspace{0.3cm}
$^2$Brunel University \\ Mathematical Sciences, UK \\ john.beasley@brunel.ac.uk \\
}
\end{center}

\begin{abstract}

In this paper we investigate the problem of order batching for picker routing. Our approach is applicable to warehouses (storage areas) arranged in the standard rectangular grid layout, so with parallel aisles and two or more cross-aisles. The motivation underlying our work is online grocery shopping in which orders may be composed of dozens of items. The approach presented directly addresses order batching, but uses a distance approximation to influence the batching of orders without directly addressing the routing problem.

We present a basic formulation based on deciding the orders to be batched together so as to optimise an objective that approximates the picker routing distance travelled.  We extend our formulation by improving the approximation for cases where we have more than one block in the warehouse. We present constraints to remove symmetry in order to lessen the computational effort required, as well as  constraints that significantly  improve the value of the linear programming relaxation of our formulation. A heuristic algorithm based on partial integer optimisation of our mathematical formulation is also presented.
Once order batching has been decided we optimally route each individual picker using a previous approach presented in the literature. 

Extensive computational results for publicly available test problems involving up to 75 orders are given for both single and multiple block warehouse configurations.

\end{abstract}


\section{Introduction}

In this section we first illustrate the order batching problem considered in this paper. We then discuss the motivation behind our work and the differences between the work presented in this paper and the work presented in our previous paper \cite{valle2017}. We also detail what we believe to be the 
contribution of this paper to the literature. 
Finally we outline the structure of the paper.

\subsection{Order batching and picker routing}

\begin{figure}[!ht]
\centering
  \includegraphics[width=0.5\textwidth]{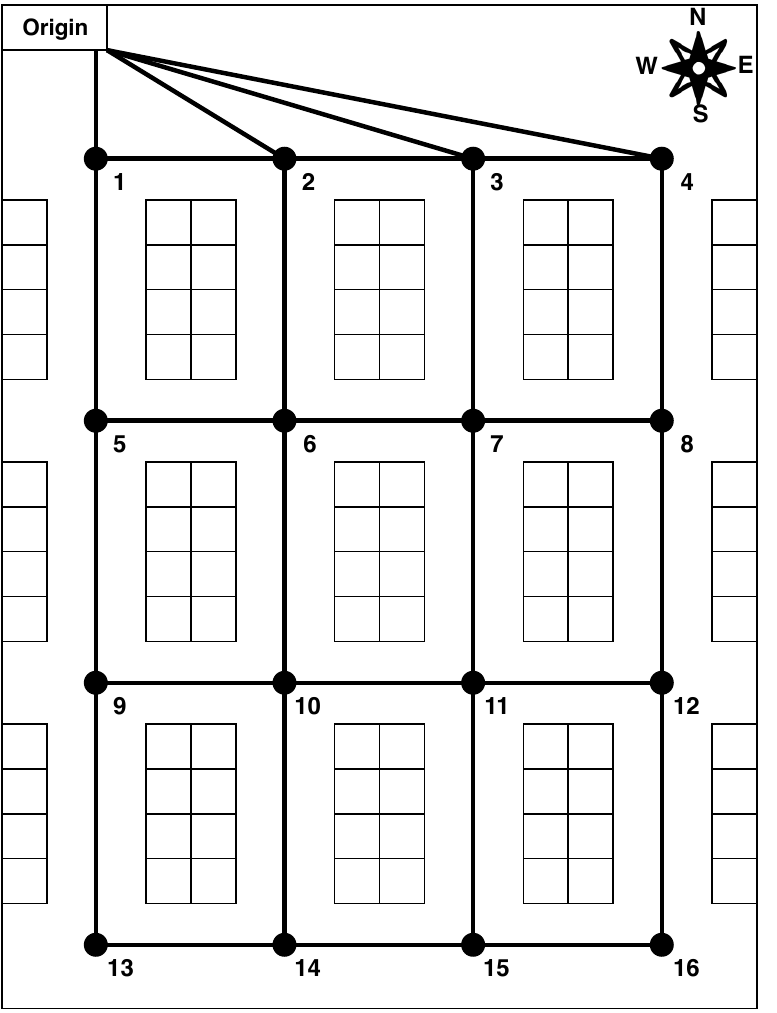}
  \caption{Warehouse structure}
  \label{fig1}
\end{figure}

To illustrate the order batching and picker routing problem consider the example warehouse shown in Figure~\ref{fig1}. 
In that figure we show a (small)  warehouse with four aisles running North$\leftrightarrow$South. 
For ease of discussion we shall use compass directions (North, South, East, West) as also shown in Figure~\ref{fig1}. 

The aisles shown are connected by cross-aisles running West$\leftrightarrow$East. So for example the first aisle contains edges $(1,5),(5,9),(9,13)$ and the first cross-aisle contains edges $(1,2),(2,3),(3,4)$. Figure~\ref{fig1} contains four aisles and four cross-aisles.
Figure~\ref{fig1} illustrates the standard warehouse rectangular grid layout that is very common,  with parallel aisles and two or more cross-aisles set at right-angles to these aisles.  

As a picker with a trolley traverses an edge in an aisle they can pick products required for the orders assigned to their trolley from the racks/shelves on either side of the aisle, these racks/shelves being shown as small squares in  
Figure~\ref{fig1}. Note here that products are only picked in aisle edges, no products are picked in 
cross-aisle edges.

In our example warehouse the origin from which empty trolleys start, and to which  trolleys filled with  products that have been picked return, is shown in the top-left hand corner of the warehouse. The warehouse structure contains a number of blocks. The first block consists of aisle edges $(1,5),(2,6),(3,7),(4,8)$; the second block consists of aisle edges $(5,9),(6,10),(7,11),(8,12)$; etc. 
Figure~\ref{fig1}  contains three blocks.
We refer to a single aisle edge as a subaisle, so for example edge $(5,9)$ is a subaisle.

Given a set of trolleys, each with known capacity, and a set of orders for products (whose locations in the warehouse where they can be picked are known) the order batching and picker routing problem is:
\begin{itemize}
\item  to batch (assign) all orders to a set of trolleys, so as to respect individual trolley capacities; and 
\item  to route each picker/trolley so as to minimise the distance travelled by that trolley as the products required for the orders so assigned are collected.
\end{itemize}
The order batching decision should be made so as to minimise the subsequent total routing distance travelled over all the trolleys used.

In terms of the routing of pickers/trolleys  two cases are distinguished in the literature: reversal and no-reversal. In reversal routing
pickers can reverse direction after traversing only part of one or more subaisles.
In no-reversal routing
pickers can only reverse direction at the end of a subaisle.

\begin{figure}[!ht]
\centering
  \includegraphics[width=0.5\textwidth]{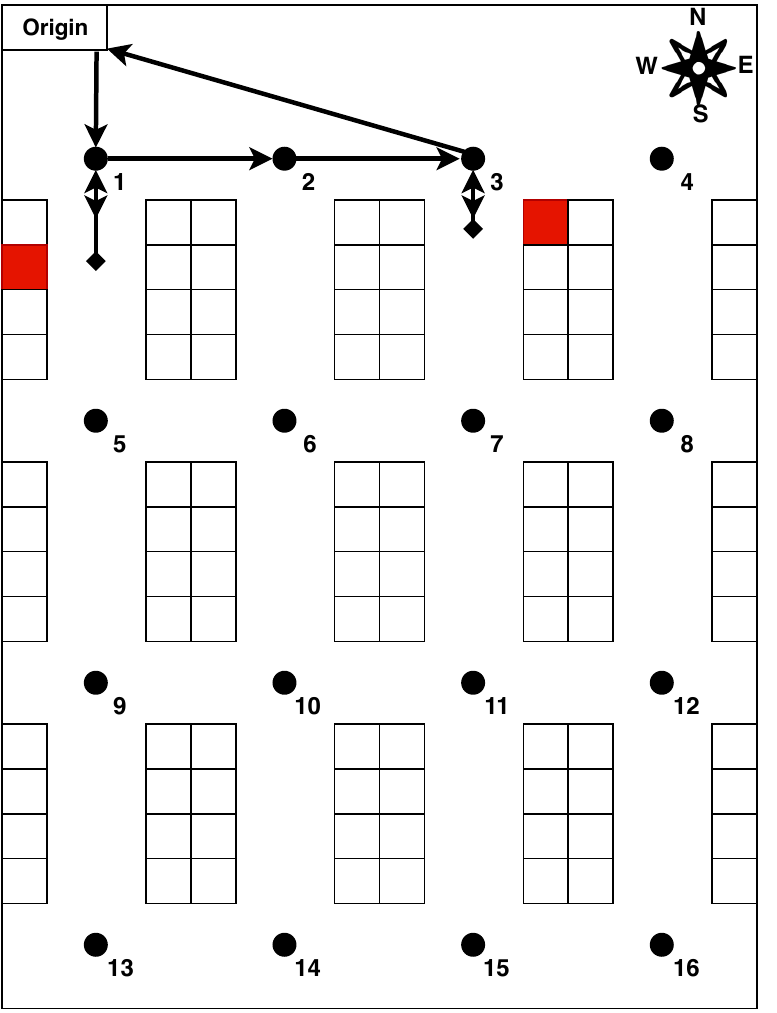}
  \caption{Trolley route: reversal}
  \label{fig2}
\end{figure}

To illustrate these two routing cases Figure~\ref{fig2} shows a case where the orders to be picked by a  single trolley only involve products situated in  subaisles $(1,5)$ and $(3,7)$, the location of these products being shown as solid squares in these subaisles. The reversal route shown involves the trolley moving down subaisle 
$(1,5)$ until it reaches the shelf where the required product is stored, it then goes back up the subaisle to vertex $1$. A similar reversal move can be seen in subaisle $(3,7)$.

\begin{figure}[!ht]
\centering
  \includegraphics[width=0.5\textwidth]{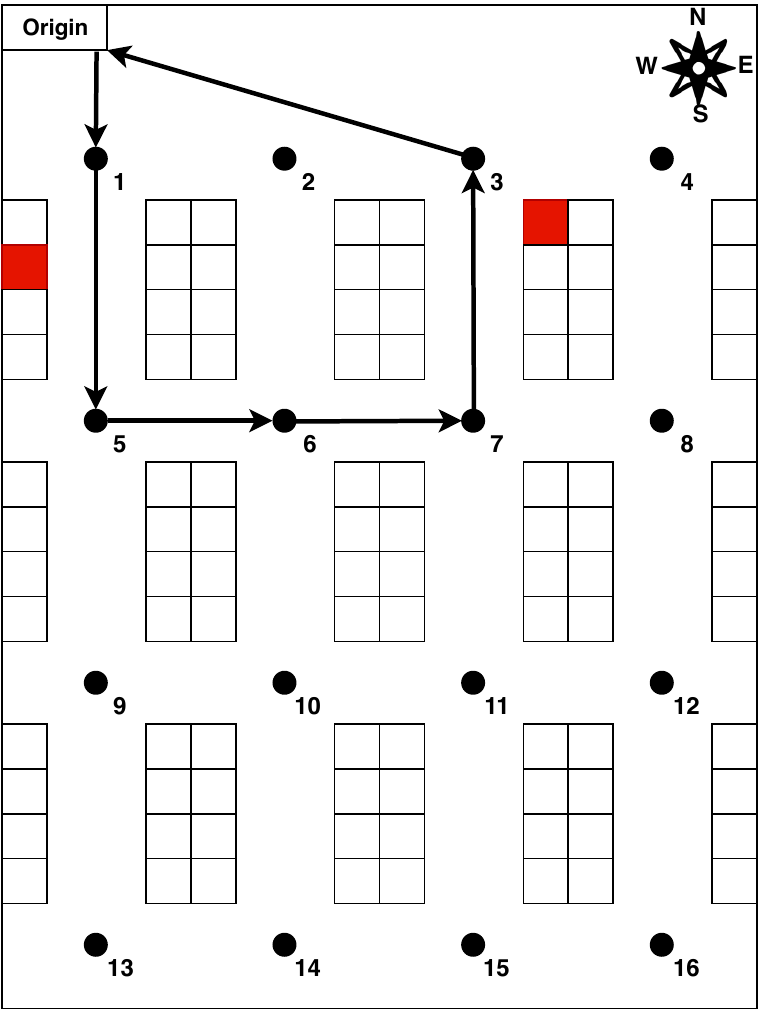}
  \caption{Trolley route: no-reversal}
  \label{fig3}
\end{figure}

Figure~\ref{fig3} shows no-reversal routing for the same situation as considered in Figure~\ref{fig2}. Here once the trolley has entered subaisle $(1,5)$ to pick the required product it must continue down the subaisle until it reaches vertex $5$,
so the cross-aisle at the end of the subaisle. Note here that for clarity of illustration we have in Figure~\ref{fig2} and Figure~\ref{fig3} arbitrarily assigned a  direction to the trolley route.

Although, as can be easily seen from Figure~\ref{fig2} and 
Figure~\ref{fig3}, allowing routes involving reversal may involve less distance than no-reversal routing it may be that
routes with no-reversal are easier to implement in practice, as they may seem more logical and intuitive for human pickers.

Previous work, \cite{valle2017}, has indicated that it is now computationally feasible to optimally route an individual trolley for relatively large problems (both for reversal and no-reversal routing).  In this paper therefore we focus on order batching. 
The approach presented directly addresses order batching, but uses a distance approximation to influence the batching of orders without directly addressing the routing problem.
We present  a  formulation that decides the orders to be batched together so as to optimise an objective that 
\emph{\textbf{approximates}} the picker routing distance travelled.  Once order batching has been decided we optimally route each individual picker/trolley using~\cite{valle2017}.

The situation considered in this paper is often referred to in the literature as a  picker-to-parts system, where a picker travels  (typically with a trolley) along aisles in a warehouse to retrieve products. As noted in~\citep{deKoster2007, marchet2015, gils2018} such systems are widely used and are believed to constitute a significant majority of all order picking systems.

\subsection{Motivation, differences and contribution}

\subsubsection{Motivation}

The motivation behind the research undertaken in this paper was \emph{\textbf{to develop an approach that could be used for the optimal batching of orders, but without explicitly routing pickers/trolleys, rather using a distance approximation to influence the batching of orders without directly addressing the routing problem.}} The key insight we had here was that because the routing of each picker/trolley occurs on a standard rectangular grid (such as seen in Figure~\ref{fig1}), then it is possible, as will become apparent below, to estimate the distance travelled to a reasonable degree of accuracy 
without explicitly deciding the picker/trolley route.

We believe that the work presented in this paper addresses a gap in the research literature. In essence, in the work surveyed later below, the quality of any particular order batching, as encountered during any solution procedure (either optimal or heuristic), is judged by constructing an \emph{\textbf{explicit}} routing for the pickers/trolleys so that they can pick the orders allocated to them. The literature appears to be almost devoid of work that evaluates the quality of any order batching in an implicit manner, without an explicit routing. The distance approximation approach presented in this paper is, however, such an implicit approach. We believe that the results presented below validate the worth of addressing this research gap.

As the literature survey given below indicates, although there have been a significant number of papers in the literature dealing with heuristics for order batching and routing, relatively little work has appeared in the literature dealing with optimal approaches.
In particular we focused on order batching since recent research, such as we presented in \cite{valle2017}, has shown that routing can (computationally) be done optimally for relatively large problems.

\subsubsection{Differences between this paper and \cite{valle2017}}

Given our previous work on the problem considered in this paper it is appropriate  to clarify here the distinct differences between the branch-and-cut approach previously presented in \cite{valle2017}
 and the distance approximation approach presented in this paper. These are:
\begin{itemize} 
\item
The approach in  \cite{valle2017} is a branch-and-cut approach where cuts (valid inequalities) are added as appropriate during the search tree. The approach presented in this paper is a 
straightforward branch and bound approach where no cuts  are added  at any stage during the search tree.

\item
The approach in  \cite{valle2017} includes cuts (valid inequalities) associated with eliminating subtours in trolley routes.  The approach presented in this paper does not involve subtours, since it does not directly address the routing problem. 

\item 
The approach in  \cite{valle2017} jointly solves the order batching and routing problem, so that it explicitly  decides simultaneously both the batching of orders  and the routes to be adopted so as to optimise over both problems (order batching and routing) directly. To achieve this it makes use of an arc based formulation.

\item  
The approach presented in this paper is an edge based formulation that directly addresses order batching, but uses a distance approximation to influence the batching of orders without directly addressing the routing problem. The advantage of this is that, as will become apparent below, problems  can be solved computationally much faster than in \cite{valle2017}.
\end{itemize}

\subsubsection{Contribution}

We believe that  the contribution of this paper to the literature is:
\begin{itemize}
\item to present a formulation of the problem of order batching based upon an approximation of the routing distance
\item to extend our formulation to provide a better distance approximation for cases where we have more than one block in the warehouse
\item to present new symmetry breaking constraints, as well as  constraints that significantly  improve the value of the linear programming relaxation of our formulation 
\item to use our formulation as a basis for a heuristic algorithm based on partial integer optimisation
\item to present results, for publicly available test problems, for both the optimal solution of our order batching formulation and for our heuristic algorithm
\item to demonstrate computationally, using the same set of test problems, that the distance approximation approach given in this paper gives results nearly as good, or better, than the approach given in our previous paper~\cite{valle2017}, but in significantly lower computation times 
\end{itemize}

\subsection{Structure of the paper}

The structure of this paper is as follows. In 
Section~\ref{sec:review} we review the literature on the order batching and picker routing problem considered. 
In Section~\ref{sec:basic} we present our basic formulation which decides the orders to be batched together so as to optimise an objective that approximates the picker routing distance travelled.
In Section~\ref{sec:extend} we extend our formulation by
improving the cross-aisle distance approximation for cases where we have more than one block in the warehouse. We also extend the aisle distance approximation and present 
constraints  to remove symmetry in order to lessen the computational effort required.
In Section~\ref{sec:practical}
we  discuss some relevant practical and computational considerations
and also present a heuristic algorithm based upon  our formulation. In Section~\ref{sec:lp} we present constraints that strengthen the linear programming relaxation of our formulation. Section~\ref{sec:results} presents our computational results for the publicly available test problems we examined. Finally in Section~\ref{sec:conclusions} we present our conclusions.

\section{Literature review}
\label{sec:review}

There is an extensive literature on the problems of order batching and picking. An early survey relating to different picking strategies was presented by~\cite{deKoster2007}, where order picking was shown to be a critical activity. In this paper we deal with batch picking, where a picker collects all of the products for one or more orders. A very comprehensive state of the art classification and review of picking systems has recently been presented by~\cite{gils2018}. 
They identify 41 papers in the literature concerned with order batching. It is clearly of little utility to repeat their work in this paper. As a consequence therefore although we, in this section, do briefly consider early work our focus is on work that is either additional to that considered 
in~\cite{gils2018}, or is 
of especial relevance to the approach we adopt.

\subsection{Early work}

Early work considered S-shaped routes~\citep{deKoster1999,  deKoster2007, hall93} where the routing rule which applies is that any aisle  in which any product is picked must be traversed in its entirety. So an aisle involving picking is either completely traversed North$\rightarrow$South from the  first cross-aisle to the  last cross-aisle, 
or completely traversed South$\rightarrow$North
from the  last cross-aisle to the  first cross-aisle.
The only exception to this rule is the last aisle visited before the trolley returns to the origin. 
Reversal is allowed in this last aisle if it aids in reducing the distance travelled.

Heuristics have been proposed for the problem of routing a single picker in warehouses with multiple blocks. An evaluation of several heuristics for single picker routing was
 given by \cite{petersen1997}. \cite{roodbergen2001_2} introduced a dynamic programming algorithm for single picker routing in warehouses with up to three cross-aisles and \cite{roodbergen2001} presented heuristics for warehouses with multiple cross-aisles. \cite{theys2010} adapted the Lin-Kernighan heuristic for the travelling salesman problem
 by \cite{lin1973} for single picker routing.

Heuristics have also been proposed for batching and routing multiple pickers. Routing distances are often estimated using single picker heuristics during the solution of the batching problem. \cite{henn2012_3} includes an extensive survey of batching methods.
For general warehouses with any number of blocks, several batching methods were presented by \cite{deKoster1999}, where routing distances are estimated using single picker methods. 


\subsection{Heuristics, single block}

\cite{lu16} presented a paper dealing with the situation where the pick-list assigned to a picker is changed (updated) during the picking operation. Computational results were presented.
\cite{chabot17} presented a paper that deals with order picking when there are weight, fragility and category constraints. They presented two formulations for the problem as well as five heuristic algorithms.  Computational results were presented for  randomly generated test problems.
\cite{hong2017} considered order batching and picker routing where pickers use S-shaped routes.
Their paper builds on their previous work~\citep{hong2012a} and involves predefining S-shaped routes. Computational results were presented for a warehouse with six aisles.

\cite{menendez17a} 
presented a variable neighbourhood search heuristic for order batching where each order has a specified due date.  They aim to batch orders and route the batches so as to minimise total tardiness.
Picker routing was done using the combined strategy 
given by~\cite{roodbergen2001}.
Computational results were presented for a large number of  test problems.
\cite{menendez17b} 
presented a two-stage variable neighbourhood search heuristic for order batching where picker routing was done using  an approach based on a modification of the combined strategy given by
~\cite{roodbergen2001}.
Computational results were presented for a large number of  test problems.
\cite{menendez17} presented both sequential and  parallel variable neighbourhood search algorithms for order batching so as to attempt to minimise the length of the longest  picker route. They categorise the problem as a min-max order batching problem, which has been previously considered in the literature by~\cite{gademann01}. Picker routing was done using both the combined strategy given by~\cite{roodbergen2001} and the dynamic programming approach of~\cite{ratliff1983}.
Computational results were presented.

\cite{chabot2018} considered an order picking problem arising in a specific industrial context which involved narrow aisles. They showed that the specific problem considered can be modelled as a vehicle routing problem. An arc reduction procedure was presented as well as two heuristic procedures and an exact branch-and-cut algorithm. Computational results were presented for randomly generated test problems.
\cite{weidinger2018} considered the situation where the products to be picked are scattered amongst the warehouse shelves. They presented a complexity proof as well as three heuristic procedures for routing a single picker. Computational results were presented for randomly generated test problems.
\cite{zulj18} presented a heuristic based on adaptive large neighbourhood search and tabu search for the order batching problem.
Pickers were routed using the  S-shaped   route strategy and the largest gap routing strategy~\citep{deKoster1999,  deKoster2007, 
hall93}. Computational results were presented for a large number of  test problems.

\cite{weidinger2019} considered the situation where the products to be picked are scattered amongst the warehouse shelves. They considered multiple locations for the warehouse (depots) where finished orders can be handed over, and new orders can be initiated. They presented a mixed-integer program for the problem of routing a single picker as well as a number of heuristic algorithms. Computational results were presented for randomly generated test problems.


\subsection{ Heuristics, multiple blocks}

 For warehouses with two cross-aisles, a variable neighbourhood search heuristic was proposed by \cite{albareda2009} and tabu search and hill climbing heuristics  by \cite{henn2012}. The latter was adapted by \cite{henn2013} to a problem that also considers order sequencing. 
\cite{kulak2012} presented an approach for joint order batching and picker routing  based on tabu search combined with a clustering approach. Order clustering (batching)  was  based on seed orders, route similarity and regret.  Computational results were presented for a number of randomly generated test problems involving both a single block and multiple blocks.

\cite{azadnia2013} considered heuristics that solve sequencing and batching first, then routing in a second stage. For warehouses with three cross-aisles, 
\cite{chen2013} focused on preventing the congestion that occurs when there are just two pickers.  For example congestion might occur if the two  pickers attempt  to simultaneously traverse the same narrow pick aisle from  opposite directions. They presented a routing heuristic 
(A-TOP) based on ant colony optimisation~\citep{dorigo1996}. Computational results were presented for their heuristic, as well as for   a modification to S-shaped routes designed to avoid congestion.

\cite{lam2014} introduced an integer programming formulation for batching, where routing distances are estimated. The problem is solved with a heuristic based on fuzzy logic.
\cite{matusiak2014} proposed a simulated annealing algorithm which includes precedence constraints with the  routing of candidate batches performed using the A* algorithm from \cite{hart1968}.
\cite{matusiak2017} considered how to incorporate skill differences between pickers in order to help minimise total order picking time. They make use of regression to estimate picker skill and use a number of heuristic procedures, based on adaptive large neighbourhood search,  to batch orders  and route pickers. Computational results were presented using data taken from a practical application.

\cite{santis2018} presented an algorithm for picker routing based upon ant colony optimisation~\citep{dorigo1996} in conjunction with a shortest path algorithm based on~\citep{floyd62, warshall62}. Computational results were presented for test problems involving up to ten aisles and four blocks.
The situation where a warehouse is divided into picking zones was considered by \cite{huang2018} who
 presented a nonlinear model for order batching that considers balancing workload between picking zones. They applied a bi-objective genetic algorithm for the solution of this model. They also presented an order batch sequencing model solved via transformation to a pseudo-boolean optimisation problem.  \cite{gils18a} presented a paper examining the relationship between storage, batching, zoning and picking using real-world data for a large warehouse in Belgium comprising 16 aisles and two blocks. They used full factorial analysis of variance to investigate a number of hypotheses.

\cite{celik2019} considered the routing of a single picker and presented complexity results. They also presented a graph-theory based merge-and-reach heuristic for the problem, with the solution obtained by this heuristic improved using the three-opt local search procedure. Computational results were presented for a number of randomly generated test problems.

\subsection{Exact approaches}

Very few exact algorithms for the joint batching and routing problem have been proposed in the literature. \cite{hong2012} dealt with multiple pickers and a single block. They considered aisle congestion where pickers are delayed by the presence of other pickers in the same aisle. In terms of picker routing they considered each picker to follow the same route around all of the aisles, with the proviso that some aisles can be skipped if desirable.  
 \cite{hong2012a} presented a formulation that depends on explicitly enumerating all possible picker routes, but no computational results for that formulation were presented, rather lower bounds on the solution to that formulation were developed as well as a heuristic procedure.
\cite{scholz2016} dealt with just a single picker and a single block. Although they indicated that their approach can be extended to more than one block, no numeric results were given.

\cite{boysen17} considered a mobile rack warehouse, where the racks on which products are stored need to be moved to create one or more aisles down which pickers can travel. In situations such as this the time to move the racks to create the aisles must be considered. They discussed both exact and heuristic approaches and presented a heuristic based on simulated annealing. Computational results were presented for  randomly generated test problems for warehouse configurations involving  just a single block. Note that, although they consider pickers and racks, their paper differs from the majority of the work considered here where the racks are fixed, not mobile.

\cite{valle2017} presented a formulation for jointly deciding order batching and picker routing 
so as to minimise the total distance travelled. Their formulation is arc based, so they explicitly distinguish between North$\rightarrow$South and South$\rightarrow$North traversals in each subaisle and West$\rightarrow$East and East$\rightarrow$West traversals in each cross-aisle. As well as introducing a number of valid inequalities based on their arc formulation they deal with subtours via an exponential number of connectivity constraints which are introduced as required in a  branch-and-cut fashion~\citep{padberg1987, padberg1991}. 
In terms of the routes adopted by pickers they considered both reversal and no-reversal routing.
Computational results were presented for a  number of publicly available test problems.

\section{Basic formulation}
\label{sec:basic}

In this section we present our basic formulation based on deciding the orders to be batched together so as to optimise an objective that approximates the picker routing distance travelled. 
In this formulation we do not attempt to construct any routing of the trolleys, rather we approximate the distance they travel. In this way we aim to get a  good assignment of orders to pickers/trolleys such that we can later optimally route each trolley individually using previous work~\citep{valle2017} presented in the literature.

In this section we also illustrate for a very small example how our formulation approximates distance without explicitly deciding trolley routes.

\subsection{Notation}
Let $W_A$ be the number of aisles and $W_B$ be the number of blocks. 
Let the origin representing the starting and return point for trolleys be labelled 0. We assume, without significant loss of generality, that this origin is located at the top left corner of the grid layout, as in Figure~\ref{fig1}.

Let the vertices in the first cross-aisle be labelled $1,2,\ldots,W_A$ in a West$\rightarrow$East direction, the vertices in the second cross-aisle are labelled $W_A+1, W_A+2,\ldots,2W_A$
in a West$\rightarrow$East direction, 
 etc. Let $E_a$ be the set of the edges in aisle $a,~a=1, \ldots, W_A$, so that $E_a=[(a+(b-1)W_A,a+bW_A)~|~b=1,\ldots,W_B]$. 
Let $E$ be the complete set of aisle edges (so $E=\bigcup_{a=1,\ldots,W_A} E_a$). Note here that we adopt an edge formulation, so we do not distinguish between a trolley traversing an aisle edge $(i,j)$ in a North$\rightarrow$South or a South$\rightarrow$North direction. Similarly we do not distinguish between a trolley traversing a 
cross-aisle edge $(i,j)$ in a West$\rightarrow$East or an East$\rightarrow$West direction.
Let the distance between any two adjacent vertices $i$ and $j$ in the grid layout be $d_{ij}$.
We shall assume throughout this paper that reasonable assumptions consistent with the underlying rectangular grid being on the Euclidean plane apply  to the distance metric, 
e.g.~non-negativity; triangle inequality satisfied; symmetric, i.e.~$d_{ij}=d_{ji}$.
We shall also assume that  the rectangular grid on which routing occurs has no missing edges (so that it is always possible to find a routing from the origin around any set of subaisle edges and then return to the origin).

Let $O$ be the set of orders that must be collected. 
Each order $o \in O$ is comprised of a set of products $P_o$ all of which have to be collected by the trolley to which that order is assigned. Let $P$ be the complete set of products involved in the complete set of orders, so that $P=\bigcup_{o \in O} P_o$. 
Let $Q(p)$ be the set of aisle edges $(i,j) \in E$ such that it is possible to collect product $p$ ($p \in P$) by traversing any of the edges in $Q(p)$. If $|Q(p)|=1$ then there is a single unique edge in $E$ where it is possible to pick product $p$. If $|Q(p)| \geq  2$ then there are two or more edges in $E$ where it is possible to pick product $p$.

Let $T$ be the number of trolleys which are available to pick the orders. 
In the formulation presented  below we are not forced to use all $T$ trolleys, and so we may use less than $T$ trolleys if it is desirable in terms of the distance travelled.
Let $B$ be the number of baskets that a trolley can carry and let $b_o$ be the known number of baskets needed to carry order $o \in O$. We assume that a basket will only contain products from a single order, even if it is partially empty, i.e.~it is not possible to put products from different orders into the same basket.

Our formulation involves a significant number of decision variables, as below, let:
\begin{itemize}
\item $z_{ot}=1$ if trolley $t$ picks order $o \in O$, zero otherwise 
\item $x_{tij}=1$ if aisle edge $(i,j)$ is traversed by trolley $t$, zero otherwise
\item $\alpha_{t}=1$ if  trolley $t$ is used to pick at least one order, zero otherwise
\item $\beta_{1ta}=1$ if the first edge in aisle $a$ in block one, so edge $(a,a+W_A)$, is traversed by trolley $t$, zero otherwise
\item $\gamma^{F1}_{ta}=1$ if trolley $t$ is used and 
aisle $a$ is the lowest indexed aisle visited by trolley $t$ in block one, zero otherwise
\item $\gamma^{L1}_{ta}=1$ if trolley $t$ is used and 
 aisle $a$ is the highest indexed aisle visited by trolley $t$  in block one, zero otherwise
\item $D^{WE}_t$ ($\geq 0$) be the contribution to the distance approximation objective function from trolley $t$ traversing cross-aisles
\end{itemize}

\subsection{Formulation}
Our basic formulation is therefore:

\begin{optprog}
\optaction[]{min} &  \objective{\sum_{t=1}^T \Bigg[ \sum_{(i,j) \in E} d_{ij} x_{tij} + \sum_{a=1}^{W_A} (d_{0a}\gamma^{F1}_{ta} + d_{a0}\gamma^{L1}_{ta})  +   D^{WE}_t \Bigg] \label{eq1} }\\

subject to: & \sum_{t=1}^T z_{ot} & = & 1 & \forall o \in O \label{eq2} \\

& z_{ot}   & \leq & \alpha_t & \forall o \in O,~t=1, \ldots, T  \label{eq3} \\

& \alpha_t & \leq & \sum_{o \in O} z_{ot} & t = 1, \ldots, T  \label{eq4} \\ 

& \sum_{o \in O} b_o z_{ot} & \leq & B \alpha_{t} & t=1, \ldots, T \label{eq5} \\

& \alpha_t & \geq & x_{tij}   & \forall (i,j) \in E,~t=1, \ldots, T \label{eq6} \\

& \alpha_t & \leq&  \sum_{(i,j) \in E} x_{tij} & t=1, \ldots, T \label{eq7} \\

& \beta_{1ta} & = & x_{ta(a+W_A)} & a=1, \ldots, W_A, t=1, \ldots, T \label{eq8} \\

& \sum_{a=1}^{W_A} \gamma^{F1}_{ta} & = & \alpha_t & t=1, \ldots, T \label{eq9} \\

& \gamma^{F1}_{ta} & \leq & \beta_{1ta} & a=1, \ldots, W_A, ~ t=1, \ldots, T \label{eq10} \\

& \gamma^{F1}_{ta} & \geq & \beta_{1ta} -\sum_{e=1,~e <a}^{W_A} \beta_{1te} & a=1, \ldots, W_A, ~ t=1, \ldots, T \label{eq11} \\

& \sum_{a=1}^{W_A} \gamma^{L1}_{ta} & = & \alpha_t & t=1, \ldots, T  \label{eq12} \\

& \gamma^{L1}_{ta} & \leq & \beta_{1ta} & a=1, \ldots, W_A, ~ t=1, \ldots, T \label{eq13} \\

& \gamma^{L1}_{ta} & \geq & \beta_{1ta} -\sum_{e=1,~e >a}^{W_A} \beta_{1te} & a=1, \ldots, W_A, ~ t=1, \ldots, T \label{eq14} \\

& \sum_{(i,j) \in Q(p)}x_{tij} & \geq & z_{ot} & \forall o \in O,~\forall p \in P_o,~t=1,\ldots,T \label{eq15} \\


& D^{WE}_t & \geq & \sum_{a=1}^{W_A} (d_{1a}\gamma^{L1}_{ta} - d_{1a}\gamma^{F1}_{ta})  & t=1, \ldots, T \label{eq16} \\
& \alpha_t & \in & \{0,1\} & t=1, \ldots, T \label{eq17} \\
& z_{ot} & \in & \{0,1\}  &  \forall o \in O,~t=1, \ldots, T \label{eq18} \\
& x_{tij} & \in & \{0,1\} & t=1, \ldots, T, ~\forall (i,j) \in E \label{eq19} \\
& \beta_{1ta}, \gamma^{F1}_{ta}, \gamma^{L1}_{ta} & \in & \{0,1\}  & t=1, \ldots, T,~a=1,\ldots,W_A \label{eq20}
\end{optprog}

In Equation~(\ref{eq1}) the term $\sum_{(i,j) \in E} d_{ij} x_{tij}$ captures the distance travelled in aisle edges by trolley $t$. The term $\sum_{a=1}^{W_A} (d_{0a}\gamma^{F1}_{ta} + d_{a0}\gamma^{L1}_{ta})$ captures the distance   associated with trolley $t$ travelling from the origin to the first cross-aisle and back. This is the distance   associated with trolley $t$ travelling from the origin to the aisle associated with a $\gamma^{F1}_{ta}$ value of one, and then travelling back to the origin from the aisle associated with a 
$\gamma^{L1}_{ta}$ value of one. Since we can arbitrarily set a direction for the trolley route (as we assume distance is symmetric) this is a valid expression. The term  $ D^{WE}_t$ captures the distance travelled in cross-aisles by trolley $t$.

Equation~(\ref{eq2}) ensures that each order is assigned to a trolley. Equation~(\ref{eq3}) ensures that an order cannot be assigned to a trolley unless that trolley is used, whilst 
Equation~(\ref{eq4}) ensures that a trolley is not used if no order is assigned to it. Equation~(\ref{eq5}) ensures that the orders assigned to any trolley $t$ cannot exceed its effective capacity $B\alpha_t$.  Equation~(\ref{eq6}) ensures that no edge $(i,j)$ can be traversed  by trolley $t$ unless that trolley is used (so $\alpha_t=1$). 
Equation~(\ref{eq7}) ensures that the trolley is not used 
if no edges are traversed by it. Note here that, strictly, Equation~(\ref{eq7}) is redundant in our formulation. However we have introduced it here to be consistent with the extension to the formulation that is presented later below.

Equation~(\ref{eq8}) sets $\beta_{1ta}$ equal to the value of the appropriate edge in block one. As above for Equation~(\ref{eq7}) this constraint could be eliminated by algebraic substitution, but we have introduced the
 $\beta_{1ta}$ variables here to be consistent with the extension to the formulation that is presented later below.

Equations~(\ref{eq9})-(\ref{eq11}) relate to the definition of the 
$\gamma^{F1}_{ta}$ variables. 
Recall here that, as defined above, the $\gamma^{F1}_{ta}$ variables are 
zero-one variables that are one if trolley $t$ is used and 
aisle $a$ is the lowest indexed aisle visited by that trolley in the first block (block one).
Equation~(\ref{eq9}) ensures that only one aisle can be picked as the first aisle (lowest indexed aisle)
in the first block if trolley $t$ is used. If the trolley is not used  then this constraint ensures that the $\gamma^{F1}_{ta}$ variables ($a=1,\ldots,W_A$) are all zero.
 Equation~(\ref{eq10}) ensures that an aisle cannot be picked as the first aisle if the first edge in the aisle is not used by trolley $t$. Equation~(\ref{eq11}) ensures that $\gamma^{F1}_{ta}$ is forced to be one if and only if the first edge in aisle $a$  in block one is used by trolley $t$ and no aisles to the West of aisle $a$ (so with index $e < a$)
 have an edge in  block one that is used.

Equations~(\ref{eq12})-(\ref{eq14}) relate to the definition of the 
$\gamma^{L1}_{ta}$ variables. 
The logic underlying these constraints is similar to that for Equations~(\ref{eq9})-(\ref{eq11}) above.
Equation~(\ref{eq12}) ensures that only one aisle can be picked as the last aisle (highest indexed aisle)
in the first block if trolley $t$ is used. If the trolley is not used  then this constraint ensures that the $\gamma^{L1}_{ta}$ variables ($a=1,\ldots,W_A$) are all zero.
  Equation~(\ref{eq13}) ensures that an aisle cannot be picked as the last aisle if the first edge in the aisle is not used by trolley $t$. Equation~(\ref{eq14}) ensures that $\gamma^{F1}_{ta}$ is forced to be one if and only if the first edge in aisle $a$ in block one is used by trolley $t$ and no aisles to the East of aisle $a$ (so with index $e > a$) have an edge in  block one that is used.

Equation~(\ref{eq15}) ensures that if order $o$ is assigned to trolley $t$, so $z_{ot}=1$, then for all products $p$ contained in order $o$ the trolley traverses at least one edge in $Q(p)$ where it is possible to collect product $p$. Note here that we do not directly decide the edge where we collect product $p$ (if there is more than one edge where the product can be collected, 
i.e.~if $|Q(p)| \geq 2$). Rather we ensure that at least one such edge is traversed. As Equation~(\ref{eq5}) ensures that we cannot assign more orders to a trolley than we have available capacity we automatically know that given an order assigned to a trolley it is feasible to pick any product $p$ in that order from any one of the edges in $Q(p)$ traversed.

Equation~(\ref{eq16}) defines 
$D^{WE}_t$, the contribution to the distance approximation objective function from trolley $t$ traversing cross-aisles. Here we approximate this by assuming that the trolley travels (in cross-aisles) from the first aisle used 
in block one
(which is given by $\gamma^{F1}_{ta}$) to the last aisle used  
in block one
 (which is given by $\gamma^{L1}_{ta}$). 
The logic underlying Equation~(\ref{eq16}) is that from Equations~(\ref{eq9}) and (\ref{eq12}), and assuming that trolley $t$ is used, only one of the $\gamma^{F1}_{ta}$ values can be one (the rest being zero) and only one of the  $\gamma^{L1}_{ta}$ values can be one (the rest being zero). So in  Equation~(\ref{eq16}) the term $\sum_{a=1}^{W_A} d_{1a}\gamma^{F1}_{ta}$
captures the cross-aisle distance from vertex one to the first aisle used in block one.  
The term $\sum_{a=1}^{W_A} d_{1a}\gamma^{L1}_{ta}$
captures the cross-aisle distance from vertex one to the last aisle used in block one. The difference therefore, 
$\sum_{a=1}^{W_A} d_{1a}\gamma^{L1}_{ta} - \sum_{a=1}^{W_A}d_{1a}\gamma^{F1}_{ta} =  
\sum_{a=1}^{W_A} (d_{1a}\gamma^{L1}_{ta} - d_{1a}\gamma^{F1}_{ta})$, captures the  cross-aisle distance between these two aisles in block one.  
If the trolley is not used then $D^{WE}_t$ will be zero. Note here that as $D^{WE}_t$ is included in the minimised objective function (Equation~\ref{eq1}) the use of an inequality in Equation~(\ref{eq16}) is valid (and adopted here because we extend the definition of $D^{WE}_t$ later below).

In Equation~(\ref{eq16})
we have assumed that the cross-aisle distances are such that  all cross-aisles distances  between any two adjacent aisles are the same. For example considering aisle 1 and aisle 2 we have
$d_{12}=d_{(W_A+1)(W_A+2)}= d_{(2W_A+1)(2W_A+2)}=,\ldots ,=d_{(W_BW_A+1)(W_BW_A+2)}$. However although we require that all  cross-aisles distances  between any two adjacent aisles are the same the entire set of cross-aisle distances need not be equal. For example the 
cross-aisle distance between aisles 1 and 2 could be different from the cross-aisle distance between aisles 2 and 3, etc.
Equations~(\ref{eq17})-(\ref{eq20}) are the integrality constraints.

\subsection{Trolley routing}

In the formulation given above we have variables related to the traversal of aisle edges. There are no variables relating to the traversal of cross-aisle edges. Since a trolley route will typically involve cross-aisle edges, as can be seen for example in Figure~\ref{fig2} and Figure~\ref{fig3}, this therefore means that in general it is \emph{\textbf{mathematically impossible}} for our formulation to 
produce a route for any trolley. Rather our formulation approximates the distance associated with a routing solution. The advantage of not needing to produce a route means that in our formulation we need not consider subtour elimination constraints (which are typically exponential in number).

To illustrate how our formulation works consider Figure~\ref{fig3}, as discussed previously above, where  the orders to be picked only involve products situated in  subaisles $(1,5)$ and $(3,7)$. Assuming, for the sake of simplicity, that all edges (aisle edges, cross-aisle edges and origin to first aisle edges) are of distance 1 then it is clear that the optimal no-reversal solution (for a single trolley) will consist of the trolley route Origin-1-5-6-7-3-Origin with the total distance travelled being 6, that route being shown in Figure~\ref{fig3}.

\begin{figure}[!ht]
\centering
  \includegraphics[width=0.5\textwidth]{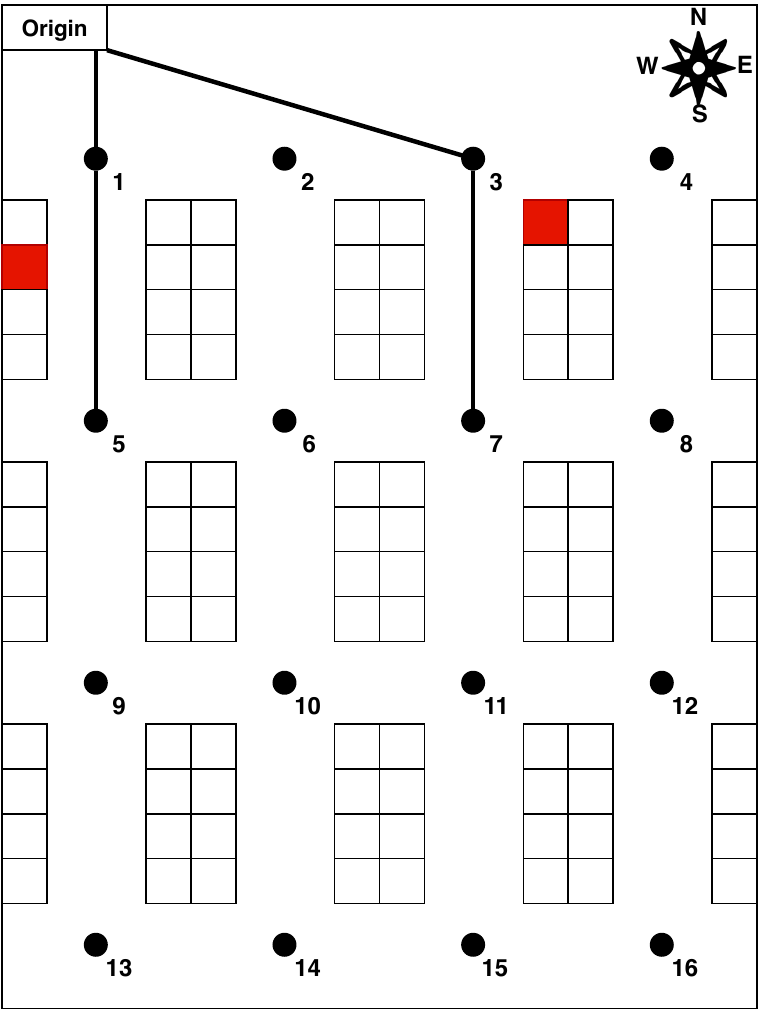}
  \caption{ Single block example: distance approximation solution}
  \label{fig4}
\end{figure}

The optimal solution to our distance approximation formulation,  Equations~(\ref{eq1})-(\ref{eq20}), will also be of value 6.
In terms of the relevant non-zero variables values that solution will have $x_{115}=x_{137}=1$, indicating that aisle edges (1,5) and (3,7) are traversed by the trolley (trolley 1).  In addition we will have $\gamma^{F1}_{11} = \gamma^{L1}_{13} =1$, indicating that the trolley has aisle 1 as the lowest indexed aisle and aisle 3 as the highest indexed aisle.  Figure~\ref{fig4} shows the distance approximation solution where the  edges shown indicate the edges involved in that solution. Note that  these edges do not constitute a trolley route.

Although the edges chosen in the distance approximation solution only correspond to a total distance of 4, the optimised distance approximation objective function, Equation~(\ref{eq1}),  will be of value 6. In detail there is a contribution of 2 from $x_{115}$ and $x_{137}$ associated with the first term in Equation~(\ref{eq1}), a contribution of 2 from $\gamma^{F1}_{11}$ and  $\gamma^{L1}_{13} $ associated with the second term in Equation~(\ref{eq1}) and a contribution of 2 from  $D^{WE}_1$ (using Equation~(\ref{eq16})) associated with the third term in Equation~(\ref{eq1}).
Note here how the  $D^{WE}_1$  term corresponds to the contribution from the cross-aisle in the optimal routing solution, 
Figure~\ref{fig3}.

So here, albeit in this simple example, we have approximated the routing distance perfectly but, as can be seen from 
Figure~\ref{fig4},
the solution to our formulation is not a valid picker/trolley route as it does not involve any cross-aisles.

\emph{\textbf{Essentially our formulation batches orders together for allocation to trolleys, as influenced by our distance approximation, but leaves the route to be adopted by each trolley to be decided later.}}

Previous work, \cite{valle2017}, has indicated that it is now computationally feasible to optimally route an individual trolley for relatively large problems.  Once order batching has been decided therefore we optimally route each individual trolley using~\cite{valle2017}, which uses an arc based formulation dealing with traversals of both aisle edges and cross-aisle edges and explicitly decides a trolley route. The optimal route (for a given trolley) is decided using a branch-and-cut approach~\citep{padberg1987, padberg1991} to eliminate subtours and ensure route connectivity. For further details as to this approach see~\cite{valle2017}.

\subsection{Discussion}

The basic formulation presented above is especially relevant when we have just one block. This is because it accurately captures in Equation~(\ref{eq16}) the West$\rightarrow$East distance in 
cross-aisles between the first aisle visited in the  block and the last aisle visited in the block. Note here that restricting attention to just one block is a situation that is commonly considered in the literature, 
e.g.~\cite{boysen17, chabot17, chabot2018, koster98, deKoster1999,
gademann01, gademann05, hong2012a, hong2017,
 lu16, 
menendez17a, menendez17b, menendez17, 
petersen1997,
rao13, weidinger2018, weidinger2019, zulj18}.

However, in this paper, we wish to consider the general situation where warehouses have more than one block. This is necessary given that the motivation underlying our work is online grocery shopping where orders may be composed of dozens of items. The picking locations for online grocery shopping almost always have more than one block (either in publicly accessible or dark stores).
It is possible to extend the formulation given above to deal with more than one block and this, along with other relevant extensions to the formulation,  is considered below.

\section{Extending the formulation}
\label{sec:extend}

The basic purpose behind much of the work presented in this section is that, although the formulation
(Equations~(\ref{eq1})-(\ref{eq20}))
presented above can be used when we have multiple blocks, we wish to  improve the quality of the distance approximation in such cases. We improve two separate components in our distance approximation objective function (Equation~(\ref{eq1})), namely:
\begin{itemize}
\item  the cross-aisle distance component ($D^{WE}_t$); and
\item the aisle distance component ($\sum_{t=1}^T  \sum_{(i,j) \in E} d_{ij} x_{tij}$).
\end{itemize}
We improve the cross-aisle component by considering the lowest and highest indexed aisles visited in any block $b \geq 2$. We improve the aisle distance component by ensuring that an even number of edges are traversed in any block.

Constraints are  also presented in this section to remove symmetry in order to lessen the computational effort required. In addition we prove that in the case of a single block 
with no-reversal routing
our distance approximation gives an optimal
allocation of orders to pickers and an optimal  routing.

\subsection{Improving the cross-aisle distance approximation}
 In the basic formulation presented above the cross-aisle distance approximation  $D^{WE}_t$ was defined in Equation~(\ref{eq16}) as being the cross-aisle distance between the first aisle used in block one and the last aisle used in block one. Now considering  Figure~\ref{fig1}  suppose that this first aisle in block one is aisle 2 and this last aisle in block one is aisle 3. It could be that a trolley, in collecting products for the orders assigned to it, also traverses an edge in aisle 1, but in some block $b \geq 2$, e.g.~edge $(9,13)$ in block 3.  This would entail additional
 cross-aisle distance in moving from aisle 2 in block one to aisle 1 in block 3 that is not currently accounted for in 
Equation~(\ref{eq16}). We can however  improve our expression for $D^{WE}_t$ to account for this in the manner indicated below. 

In the general case we are seeking the lowest indexed aisle used in any block $b \geq 2$. If that is to the West  of the lowest indexed aisle used in block one then we have an additional cross-aisle contribution to add to 
$D^{WE}_t$. A complication here however is that a trolley may only visit edges in block one and never visit any edges in a higher block, and this needs to be taken into consideration.

Define $G_1$ as the set of aisle edges in block one, so $G_1=[(a,a+W_A)~|~a=1,\ldots,W_A]$. Let $\hat{\alpha_t}=1$ if trolley $t$ visits some block $b \geq 2$, zero otherwise. Then the constraints relating this variable to the variables that we had previously above are:
\begin{optprog}
& \hat{\alpha_t} & \leq & \alpha_t  & t=1, \ldots, T \label{eq21} \\
& \hat{\alpha_t} & \geq & x_{tij}     & \forall (i,j) \in E \setminus G_1, \;\;\;  ~t=1, \ldots, T \label{eq22}  \\
& \hat{\alpha_t} & \leq & \sum_{(i,j) \in E \setminus G_1} x_{tij} & t=1, \ldots, T \label{eq23} 
\end{optprog}

Equation~(\ref{eq21}) ensures that $\hat{\alpha_t}$ is zero if trolley $t$ is not used at all. Equation~(\ref{eq22}) ensures that $\hat{\alpha_t}$ is one if any edge in any block $b \geq 2$ is used.  Equation~(\ref{eq23}) ensures that $\hat{\alpha_t}$ is zero if no  edge in any block $b \geq 2$ is used. Equations~(\ref{eq22}) and (\ref{eq23}) are equivalent to Equations~(\ref{eq6}) and (\ref{eq7}), but particularised to blocks  $b \geq 2$.

Let $\beta_{2ta}=1$ if aisle $a$ is visited by trolley $t$ in one or more of the blocks $b \geq 2$, zero otherwise. Let $\gamma^{F2}_{ta}=1$ if aisle $a$ is the lowest indexed aisle visited by trolley $t$ when considering all blocks $b \geq 2$, zero otherwise. Then the constraints relating to these variables are:
\begin{optprog}
 & \beta_{2ta} & \leq & \hat{\alpha_t} & a=1, \ldots, W_A,~ t=1, \ldots, T \label{eq24} \\
 & \beta_{2ta} & \geq & x_{tij} & \forall (i,j) \in E_a \setminus (a,a+W_A) \;\;\; a=1, \ldots, W_A,~ t=1, \ldots, T \label{eq25} \\
 & \sum_{a=1}^{W_A} \gamma^{F2}_{ta} & = & \hat{\alpha_t} &  t=1, \ldots, T \label{eq26} \\
 & \gamma^{F2}_{ta} & \leq & \beta_{2ta} & a=1, \ldots, W_A, ~ t=1, \ldots, T \label{eq27} \\
 & \gamma^{F2}_{ta} & \geq & \beta_{2ta} -\sum_{e=1,~e <a}^{W_A} \beta_{2te} & a=1, \ldots, W_A, ~ t=1, \ldots, T \label{eq28}
\end{optprog}

Equation~(\ref{eq24}) ensures that the $\beta_{2ta}$ variables are all zero if the trolley is not used in blocks $b \geq 2$. Equation~(\ref{eq25}) ensures that $\beta_{2ta}$ is set to one if any edge in aisle $a$ and some block $b \geq 2$ is traversed by trolley $t$. Equations~(\ref{eq26})-(\ref{eq28}) are equivalent to Equations~(\ref{eq9})-(\ref{eq11}) previously presented and explained above, but particularised to blocks $b \geq 2$.

Now if the lowest  indexed aisle visited in blocks $b \geq 2$ is West of the lowest index aisle visited in block one we automatically know that we can add an additional component to the cross-aisle distance approximation comprising twice the cross-aisle distance between these two aisles. For example if in Figure~\ref{fig1} the lowest indexed aisle visited in block one is aisle 2, but the lowest indexed aisle visited in any block 
$b \geq 2$ 
is aisle 1 then from aisle 2 in block one we need to use cross-aisles to get to aisle 1, but then come back again to aisle 2 in order to be able to reach the highest indexed aisle in block one.
Introduce $z^F_t \geq 0$ as the additional contribution to the distance approximation from this extra cross-aisle  distance for trolley $t$. Then we have:
\begin{equation}
z^F_t \geq 2\sum_{a=1}^{W_A} (d_{1a}\gamma^{F1}_{ta} - d_{1a}\gamma^{F2}_{ta}) 
-2d_{1W_A}(1-\hat{\alpha_t})
\;\;\;  t=1, \ldots, T 
\label{eq29}
\end{equation}

In Equation~(\ref{eq29})  the summation term represents the cross-aisle distance between the  lowest indexed aisles, namely the aisles  associated with $\gamma^{F1}_{ta} $ and  $\gamma^{F2}_{ta}$ being one. The difference $ \sum_{a=1}^{W_A} (d_{1a}\gamma^{F1}_{ta} - d_{1a}\gamma^{F2}_{ta})$ will only be positive  when either the lowest  indexed aisle visited in blocks $b \geq 2$ is West of the lowest index aisle visited in block one, or the trolley only visits block one and never visits any block $b \geq 2$. For this reason we include the term $-2d_{1W_A}(1-\hat{\alpha_t})$ to ensure that if the trolley never visits any block $b \geq 2$ (so $\hat{\alpha_t}=0$) the right-hand side of Equation~(\ref{eq29}) is $\leq  0$.  Noting that we have a minimisation objective it is therefore valid to add $z^F_t$ to the right-hand side of Equation~(\ref{eq16})  which defines the 
cross-aisle contribution $D^{WE}_t$ to the distance approximation associated with trolley $t$.

Above we have discussed the case where the lowest  indexed aisle visited in blocks $b \geq 2$ is West of the lowest index aisle visited in block one. However an analogous situation occurs when the highest indexed aisle visited in blocks $b \geq 2$ is East of the highest indexed aisle visited in block one, so in this case we can also add an additional contribution to the cross-aisle distance approximation.

Proceeding in a very similar fashion as above  let $\gamma^{L2}_{ta}=1$ if aisle $a$ is the  highest indexed aisle visited by trolley $t$ when considering all blocks $b \geq 2$, zero otherwise. Then we have the constraints:
\begin{optprog}
& \sum_{a=1}^{W_A} \gamma^{L2}_{ta} & = & \hat{\alpha_t} &  t=1, \ldots, T \label{eq30} \\
& \gamma^{L2}_{ta} & \leq & \beta_{2ta} & a=1, \ldots, W_A, ~ t=1, \ldots, T \label{eq31} \\
& \gamma^{L2}_{ta} & \geq & \beta_{2ta} -\sum_{e=1,~e >a}^{W_A} \beta_{2te} &  a=1, \ldots, W_A, ~ t=1, \ldots, T \label{eq32}
\end{optprog}

Equations~(\ref{eq30})-(\ref{eq32})  are equivalent  to Equations~(\ref{eq12})-(\ref{eq14}) previously presented and explained above, but particularised to blocks $b \geq 2$. Introducing $z^L_t \geq 0$ as the further additional contribution to the distance approximation from this extra cross-aisle  distance for trolley $t$ we have:
\begin{equation}
z^L_t \geq 2\sum_{a=1}^{W_A} (d_{1a}\gamma^{L2}_{ta} - d_{1a}\gamma^{L1}_{ta}) 
 \;\;\;  t=1, \ldots, T
\label{eq33}
\end{equation}

We can therefore replace Equation~(\ref{eq16}) by
\begin{equation}
D^{WE}_t  \geq \sum_{a=1}^{W_A} (d_{1a}\gamma^{L1}_{ta} - d_{1a}\gamma^{F1}_{ta})  + z^F_t + z^L_t
 \;\;\; t=1, \ldots, T 
\label{eq34}
\end{equation}

\subsection{Improving the aisle distance approximation}

In our formulation we have a single zero-one variable $x_{tij}$ associated with each aisle edge. As such we do not distinguish between the edge being traversed in a North$\rightarrow$South direction or in a South$\rightarrow$North direction. However it is clear that for each block the number of edges traversed  
North$\rightarrow$South must equal the number of edges traversed South$\rightarrow$North (in order for the trolley to return to the origin).
If we have any edge which is traversed by the trolley in both directions, i.e.~North$\rightarrow$South and South$\rightarrow$North, then the contribution to the distance approximation made by the term $d_{ij} x_{tij}$ in Equation~(\ref{eq1}) will be an underestimate, since we will only be counting the edge distance once and not twice.

Suppose that we have an odd number of edges traversed in a block.
Then  we know that each edge cannot be traversed just once, as that would not enable the trolley to return to the origin. To ensure that we traverse an even number of edges then there may be an additional traversal. For example, either one of the chosen edges in the block could be traversed twice, or one of the currently unused edges in the block traversed once. Alternatively to ensure that we traverse an even number of edges there may be one less traversal of some currently traversed edge.

Define $G_b$ as the set of aisle edges in block $b$, so $G_b=[(a+(b-1)W_A,a+bW_A)~|~a=1,\ldots,W_A]$. Let $y_{tb}=1$ if we have an odd number of edges traversed in block $b$ by trolley $t$, zero otherwise. Let $D_{b}$ be the minimum length of any subaisle in block $b$, so $D_b= \mbox{min}[d_{ij}~|~(i,j) \in G_b]$. 

Now the number of edges involved in  the solution associated with trolley $t$ in block $b$ is
$\sum_{(i, j) \in G_b} x_{tij}$. Introduce general integer variables $w_{tb}$ where $0 \leq w_{tb} \leq \lfloor W_A/2 \rfloor$, then we have:
\begin{equation}
\sum_{(i, j) \in G_b} x_{tij} = 2w_{tb} + y_{tb} \;\;\;
t=1, \ldots ,T,~b=1, \ldots, W_B
\label{eq35}
\end{equation}

Equation~(\ref{eq35}) ensures that if the number of edges traversed in block $b$ by trolley $t$  is odd then $y_{tb}$ will be one.
We can therefore add the term $ \sum_{b = 1}^{W_B} D_{b} y_{tb}$  for trolley $t$ to the objective function, Equation~(\ref{eq1}), to give the new objective as:
\begin{equation}
\min \sum_{t=1}^T \Bigg[ \sum_{(i,j) \in E} d_{ij} x_{tij} + \sum_{a=1}^{W_A} (d_{0a}\gamma^{F1}_{ta} + d_{a0}\gamma^{L1}_{ta})  +   D^{WE}_t + 
\sum_{b = 1}^{W_B} D_{b} y_{tb} \Bigg]
\label{eq36}
\end{equation}

Here as $y_{tb}$ is involved in the minimisation objective it will be zero unless the  number of edges traversing the block is odd in which case the only way to satisfy the equality constraint, Equation~(\ref{eq35}), is to set $y_{tb}$ to one given that the $w_{tb}$ variables are general integer variables.

\subsection{Symmetry}

In our formulation we have $T$ trolleys, each with identical capacity. A consequence of this is that there are $T!$ solutions of equal (optimal) value, consisting of exactly the same sets of orders for each trolley, but where we simply permute the trolleys assigned to each such set. 

In order to help eliminate symmetry, and hence potentially lessen the computational effort required to solve the problem, we can impose the following constraints:
\begin{equation}
\sum_{t = 1}^{o} z_{ot} = 1 \;\;\; o = 1, \ldots, T
\label{eq37}
\end{equation}
\begin{equation}
z_{ot} = 0 \;\;\; o = 1, \ldots, (T-1),~~t=(o+1), \ldots ,T
\label{eq38}
\end{equation}

Equation~(\ref{eq37}) enforces the condition that the first order must be assigned to the first trolley, the second order either to trolley one or to trolley two, the third order to one of the first three trolleys, etc. Equation~(\ref{eq38}) follows as a logical consequence of Equation~(\ref{eq37}).

Although symmetry breaking constraints of the form shown in Equations~(\ref{eq37}) and (\ref{eq38}) have been seen previously in the literature  it is possible to generate stronger symmetry breaking constraints.

Considering each order $o$ in turn we can break symmetry by seeking a solution such that order $o$ is assigned to trolley $t$ if order $o$ has not been assigned to a lower indexed trolley, i.e.~if $\sum_{r=1}^{t-1}z_{or} = 0$ and trolley $t$ is completely free to assign order $o$ to as none of the previous orders $1,2,\ldots,o-1$ have  been assigned to trolley $t$, i.e.~if $\sum_{q=1}^{o-1}z_{qt} = 0$. The constraint that breaks symmetry in this fashion is:
\begin{equation}
z_{ot} \geq 1 - \sum_{r=1}^{t-1}z_{or} - \sum_{q=1}^{o-1}z_{qt}
\;\;\;o = 2, \ldots, |O|,~t=2, \ldots ,T
\label{eq39}
\end{equation}

To illustrate Equation~(\ref{eq39}) suppose $o=4$ and $t=3$. Then the constraint becomes $z_{43} \geq 1 -  z_{41} - z_{42} - z_{13}  - 
z_{23} - z_{33}$. In other words order 4 must be assigned to trolley 3 if order 4 has not already been assigned to trolleys 1 or 2 and none of orders 1,2,3 have been assigned to trolley 3. 

The logic behind this symmetry breaking constraint is that it is
possible to relabel the trolleys in the optimal allocation of orders to trolleys to satisfy Equation~(\ref{eq39}). This can be done as follows. Label the trolley to which order 1 is assigned trolley 1. Now consider the remaining trolleys and label the trolley which contains the lowest indexed order trolley 2. Now consider the remaining trolleys and label the trolley which contains the lowest indexed order trolley 3, etc. 

This relabelling will satisfy Equation~(\ref{eq39}) since trolley $t$ in this relabelling will contain an order $o$ such that no lower indexed orders $1,2,\ldots,o-1$ are in trolley $t$ and order $o$ has not been assigned to a lower indexed trolley.

As far as we are aware symmetry breaking constraint 
Equation~(\ref{eq39}) has not appeared previously in the literature. 

Clearly to make best computational use of the symmetry breaking constraints presented above we need to index the set of orders $O$ appropriately. In the computational results reported later below we indexed the orders by sorting them into decreasing $b_o$ order, ties broken by sorting them into decreasing $|\bigcup_{p \in P_o} Q(p)|$ order. So here the orders are indexed in decreasing order of the number of baskets required, ties broken by decreasing order of the number of edges involved in the order.

\subsection{Multiple block example}

\begin{figure}[!htp]
\centering
  \includegraphics[width=0.5\textwidth]{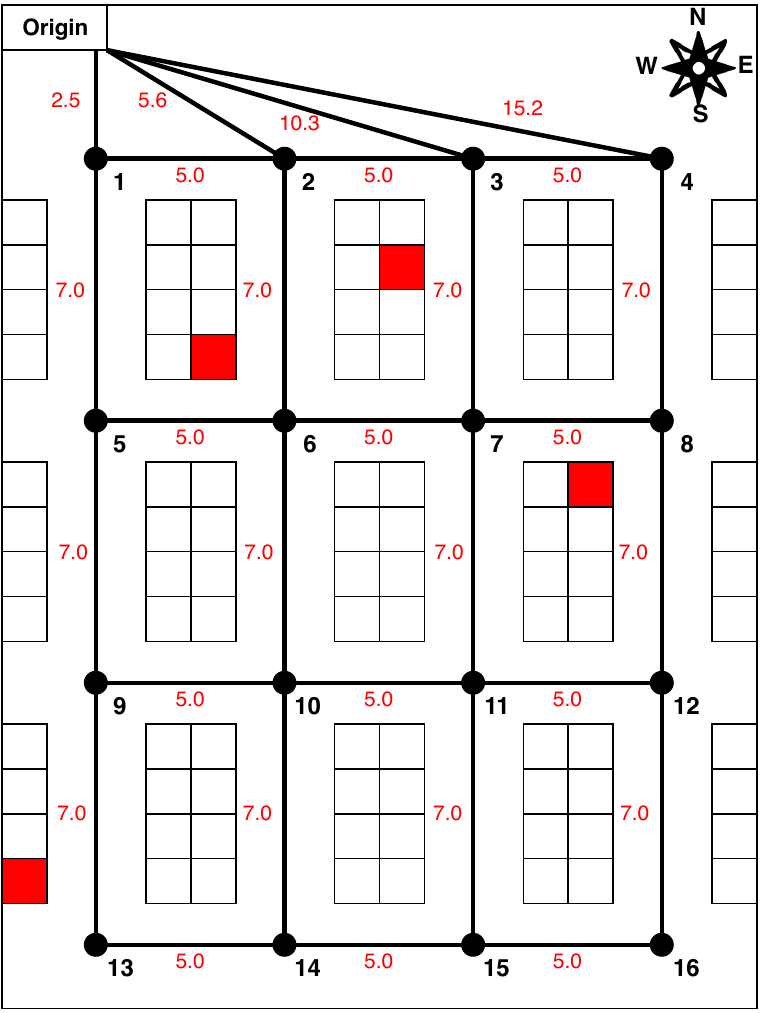}
  \caption{Warehouse: multiple block example}
  \label{fig5}
\end{figure}

To illustrate the application to multiple blocks of the distance approximation extension presented in this section Figure~\ref{fig5} shows a single picker/trolley multiple block example (with three blocks), using the same notation as in Figure~\ref{fig2}, but now with distances shown next to every edge. In Figure~\ref{fig5} the orders to be picked by a  single trolley involve products situated in  subaisles $(2,6)$, $(3,7)$, $(8,12)$ and $(9,13)$, the location of these products being shown as solid squares in these subaisles. 

\begin{figure}[!htp]
\centering
  \includegraphics[width=0.5\textwidth]{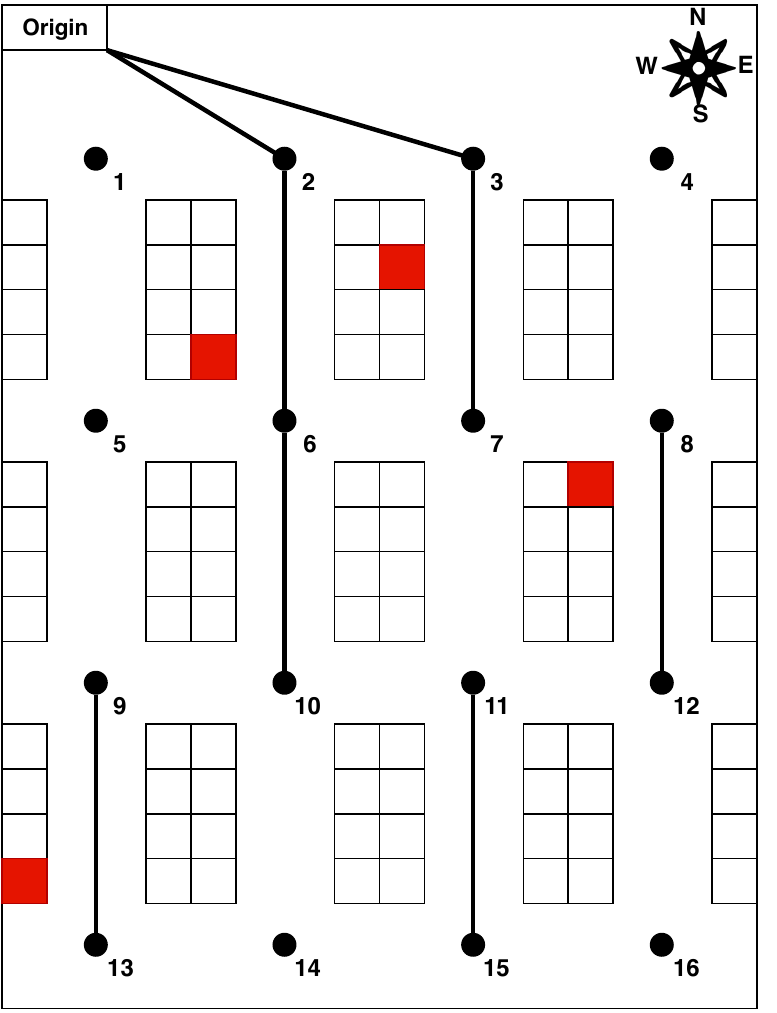}
  \caption{Multiple block example: distance approximation solution}
  \label{fig6}
\end{figure}

\begin{figure}[!htp]
\centering
  \includegraphics[width=0.5\textwidth]{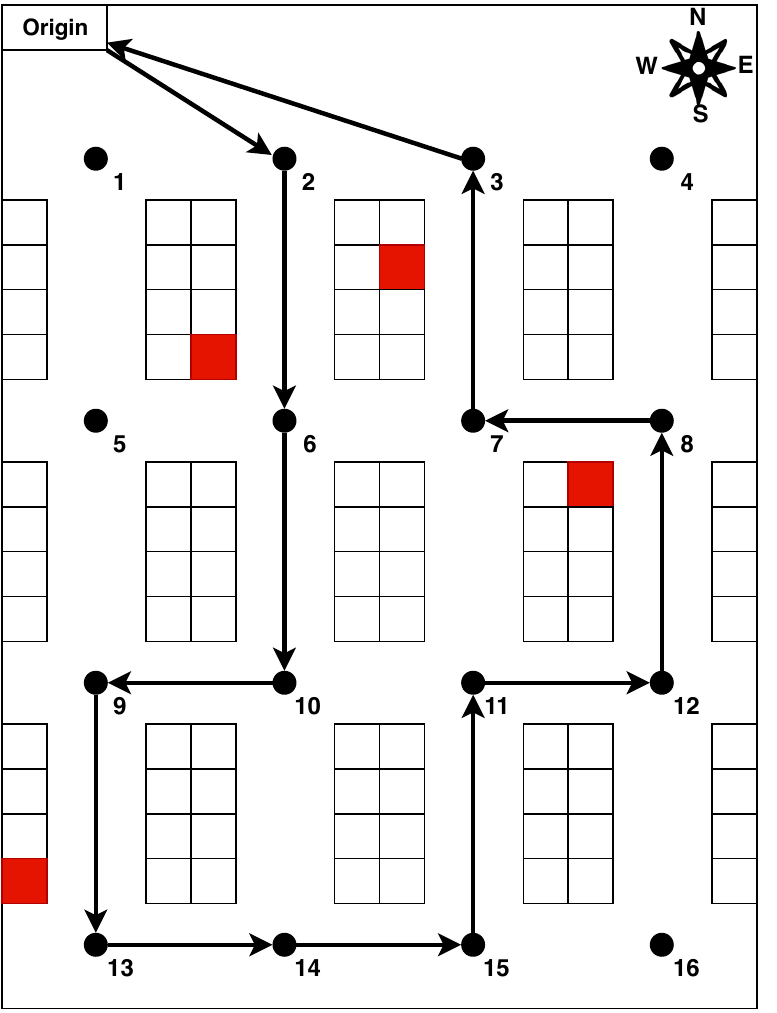}
  \caption{ Multiple block example: optimal no-reversal solution}
  \label{fig7}
\end{figure}

Figure~\ref{fig6} shows the distance approximation solution, note  that as in Figure~\ref{fig4} this is not a picker/trolley route. Figure~\ref{fig7} shows the optimal no-reversal route for this example and by reference to Figure~\ref{fig5} this involves a total distance of 82.9.  The value given by our distance approximation is also 82.9, so here we have approximated the routing distance perfectly, but without explicitly routing the picker/trolley. In detail, and referring to Figure~\ref{fig5}, Figure~\ref{fig6} and 
Equations~(\ref{eq34})-(\ref{eq36}), there is:
\begin{itemize}
 \item a contribution of 28.0 from the subaisle edges from which we pick products (the first term in Equation~(\ref{eq36})), relating to edges (2,6), (3,7), (8,12) and (9,13) in Figure~\ref{fig6},
 \item a contribution of 15.9 from the origin to first aisle vertices (the second term in Equation~(\ref{eq36})) as in the optimal solution $\gamma^{F1}_{12} = \gamma^{L1}_{13} = 1$, represented in 
Figure~\ref{fig6} as the edges from the origin to vertices 2 and 3 respectively,
 \item a contribution of 25.0 from the $D^{WE}_1$ term (the third term in Equation~(\ref{eq36})), representing the cross-aisle contribution, which can be further split (see 
Equation~(\ref{eq34})), into a contribution of 5.0 in moving from aisle 2 to aisle 3; 10.0 in moving from aisle 2 to aisle 1 and back ($z_1^F = 10.0$); and  10.0 in moving from aisle 3 to aisle 4 and back ($z_1^L = 10.0$); not represented in Figure~\ref{fig6},
 \item  a contribution of 14.0 from the fourth term in Equation~(\ref{eq36}) as  $y_{12} = y_{13} = 1$ so that an even number of subaisle edges are traversed in the second and third blocks (enforced by Equation~(\ref{eq35})). These are represented, for illustration purposes, as edges (6,10) and (11,15) in Figure~\ref{fig6}. 
\item this gives a total of $28.0+15.9+25.0+14.0 = 82.9$. 
\end{itemize}

Recall here that, as mentioned above, in our distance approximation approach we batch orders together so as to approximate the picker routing distance travelled.  Once order batching has been decided we optimally route each individual picker/trolley using~\cite{valle2017}. The question therefore arises as to the connection between the structure of the distance approximation solution associated with any individual picker (which, whilst not a route, does involve some edges/subaisles) and the  subsequent 
no-reversal route for that picker.

We can illustrate this connection using Figure~\ref{fig6} and Figure~\ref{fig7}. The no-reversal picker route adopted as associated with the distance approximation solution shown in Figure~\ref{fig6} must contain any subaisle edge associated with any picked product $p$ that can only be picked in one subaisle, i.e.~$|Q(p)|=1$. As all of the test problems used in this paper have $|Q(p)|=1~\forall p \in P$ the picker route associated with the distance approximation solution shown in 
Figure~\ref{fig6} must contain the subaisles
$(2,6)$, $(3,7)$, $(8,12)$ and $(9,13)$. 

It will be seen that the no-reversal picker route in 
Figure~\ref{fig7} does contain these subaisles. All other edges associated with the distance approximation solution may, or may not, be included in the picker route. In part this is because there may be alternative routing solutions (involving different edges) of equal distance. For example referring to   Figure~\ref{fig5} and Figure~\ref{fig7} it is clear that a number of alternative routing solutions of equal value not involving subaisle $(6,10)$ and/or subaisle $(11,15)$  exist, 
e.g.~with regard to subaisle $(6,10)$ there is an alternative  routing solution of equal value  where edges $(6,10)$ and $(10,9)$ in
Figure~\ref{fig7}
are replaced by edges $(6,5)$ and $(5,9)$.

\subsection{Distance approximation clarification}

Our distance approximation formulation is to optimise~Equation~(\ref{eq36}) subject to 
Equations~(\ref{eq2})-(\ref{eq15}),(\ref{eq17})-(\ref{eq35}),(\ref{eq37})-(\ref{eq39}). For any particular problem instance solving this formulation to proven optimality yields an optimal distance approximation objective function value ($V_1$, say)  and too provides an allocation of orders to pickers/trolleys. Each trolley is then individually optimally routed using \cite{valle2017} giving a total routing distance ($V_2$, say). 

Note how in our approach we separate the allocation of orders to trolleys (so order batching) from the subsequent routing. Suppose, for the instance under consideration, we were also to jointly decide order batching and picker routing so as to minimise the total distance travelled (that distance being $V_3$, say). Since $V_3$ pertains to the situation where we
 simultaneously decide  both the batching of orders  and the routes to be adopted so as to optimise over both problems (order batching and routing) directly
 it corresponds to the minimal routing distance that can be achieved.

In order to clarify the relationship between these three distances, $V_1$,  $V_2$ and  $V_3$, we have that in the case of no-reversal routing:
\begin{itemize}
\item for a single block instance $V_1 = V_2 = V_3$. In other words our distance approximation gives \emph{\textbf{an optimal allocation of orders to trolleys and an optimal routing}}. 
\item for a multiple block instance $V_1 \leq V_2$ and $V_2 \geq V_3$.  In other words our distance approximation solution $V_1$ provides a \emph{\textbf{lower bound on the subsequent routing distance $V_2$ which is itself an 
upper bound  on the minimal routing distance $V_3$}}.
\end{itemize}

Although our distance approximation formulation is based upon assuming that routing is done in a no-reversal fashion (as illustrated in  Figure~\ref{fig3}) it is the case that, since it provides an allocation of orders to trolleys, these trolleys can be routed in a reversal fashion (as illustrated in Figure~\ref{fig2}). In the case of reversal routing  we have:
\begin{itemize}
\item $V_2  \geq V_3$. In other words our distance approximation 
gives an allocation of orders to trolleys whose routing provides an \emph{\textbf{upper bound $V_2$ on the minimal routing distance $V_3$}}.
\end{itemize}

Proving that for a single block instance with 
no-reversal routing $V_1 = V_2$, is simply done.  Recall that our distance metric is consistent with the underlying rectangular grid (Figure~\ref{fig1})  being on the Euclidean plane. Suppose our distance approximation solution involves an even number of subaisle edges, for example as in Figure~\ref{fig4}. Then a route with the same value as the approximation solution is given by travelling from the origin to the West-most subaisle, traversing that subaisle in a  
North$\rightarrow$South direction, moving to the next subaisle and traversing that subaisle in a 
South$\rightarrow$North direction, moving to the next subaisle and traversing that subaisle in a 
North$\rightarrow$South direction, etc; finally returning to the origin from the Northern end of the last subaisle traversed. If the distance approximation solution involves an odd number of subaisle edges then we find a route as for an even number of edges, but simply re-traverse the last subaisle visited in a 
South$\rightarrow$North direction and then return to the origin. Recall here that if we have an odd number of edges then the distance approximation has terms to account for this (see Equations~(\ref{eq35}),(\ref{eq36})).
Hence we have $V_1 = V_2$.
In such cases we must also have $V_2=V_3$, i.e.~we have an optimal allocation of orders to trolleys, since if not then our distance approximation solution would not be optimal.

\section{Practical and computational considerations and a heuristic}
\label{sec:practical}

There are a number of practical and computational issues concerned with our formulation that are worthwhile highlighting here. These are considered below.
We also present a heuristic algorithm directly based upon  our mathematical distance approximation formulation.

\subsection{Practical considerations}
In the formulation presented above we, in common with the vast majority of the published literature, focus on batching orders so as to minimise the total distance travelled. Clearly any formulation may need amendment for specific conditions arising in practical situations. In terms of the problem considered in this paper two such conditions might be:
\begin{itemize}
\item the height of the storage location from which the products are to be picked may affect the time required, and this needs to be taken into account 
\item it might be important to balance the workload for different trolleys
\end{itemize}

With respect to the first of these conditions the time required to pick the orders allocated to a trolley essentially becomes significant if there is a time constraint with regard to the trolley, i.e.~if there is a constraint upon the total time that a trolley can spend away from the origin. Such a constraint might arise, for example, if orders have to be picked quickly and/or it is near the end of the working day. Amending our formulation to incorporate this constraint can be achieved by defining a standard trolley speed (thereby enabling a conversion from distance travelled to time taken) and incorporating a term corresponding to the time required to pick each order, thereby giving the time taken by each trolley, which can be suitably constrained. Here, because of space considerations we have left the precise mathematical changes needed to our formulation to the reader.

With respect to the second of these conditions if we define the workload associated with  a trolley $t$ as given by $\sum_{o \in O} b_o z_{ot}$, so the number of baskets in the trolley, then to balance the workload over the $T$ trolleys we might wish to minimise (max[$\sum_{o \in O} b_o z_{ot}~|~t=1,\ldots,T$] - 
min[$\sum_{o \in O} b_o z_{ot}~|~t=1,\ldots,T$]). This can be linearised in a standard way, introduce variables $\theta^{max}$ and $\theta^{min}$ defined using $\theta^{max} \geq \sum_{o \in O} b_o z_{ot}~~t=1,\ldots,T$ and $\theta^{min} \leq 
\sum_{o \in O} b_o z_{ot}~t=1,\ldots,T$ and minimise $(\theta^{max} - \theta^{min})$. In this minimisation $\theta^{max}$ is the highest workload associated with any trolley and $\theta^{min}$ is the lowest workload associated with any trolley. The objective adopted minimises the difference between these two workloads, hence achieving workload balance. If workload were to be defined based on the time taken by each trolley (or the distance travelled by each trolley) then this can be easily dealt with in a similar fashion.

\subsection{Redundant constraints}

The number of constraints  associated with Equation~(\ref{eq15}) can be reduced by eliminating redundant constraints. For each $z_{ot}$ term on the right-hand side of that constraint it is possible that some of the sets of edges considered on the left-hand side are dominated. 

For example, with reference to Figure~\ref{fig1}, suppose that we have  two products in some order $o$  such that the first of these products can be picked from edge $(2,6)$, i.e.~picked from just one subaisle, and the other product can be picked from edges $(2,6),(7,11)$, i.e.~picked from two subaisles. Then it is clear that if this order is assigned to trolley $t$ it must traverse edge $(2,6)$ to collect the first product, in which case that edge can also supply the second product and as a consequence the constraint relating to collection of the second product in Equation~(\ref{eq15}) involving edges $(2,6),(7,11)$ is redundant.

Eliminating redundant constraints for order $o \in O$ is trivial using the following algorithm. Iteratively consider each $p \in P_o$ in turn and if there exists another product $q \in P_o$ ($q \neq p$) 
such that $Q(q) \subseteq Q(p)$ then remove product $p$ from $P_o$. This removal is valid as 
both products $p$ and $q$ are in the same order and 
product $q$ has a set of edges which  are a subset of the set of edges associated with product $p$. Hence the constraint (Equation~(\ref{eq15})) relating to product $p$ is redundant as a trolley visiting an edge to pick product $q$ can automatically pick product $p$ using the same edge.

\subsection{Zero-one variables}

Our formulation contains a significant number of zero-one variables. However detailed consideration of the constraints involved reveals that provided the variables  $z_{ot}$ and $x_{tij}$ are declared as zero-one, and the variables $w_{tb}$ are declared as general integer variables, then 
 all  other zero-one variables 
will be naturally integer in the optimal solution to the problem if they are  declared as continuous variables lying between zero and one. 
For the avoidance of doubt the variables currently specified as zero-one in the discussion above that can be taken as continuous variables lying between zero and one are:
$\alpha_{t}$, $\beta_{1ta}$, $\gamma^{F1}_{ta}$, $\gamma^{L1}_{ta}$, $\hat{\alpha_t}$, $\beta_{2ta}$, $\gamma^{F2}_{ta}$, $\gamma^{L2}_{ta}$ and  $y_{tb}$.

Essentially  we are saying here that the key zero-one variables are related to the assignment of orders to trolleys and the assignment of edges to trolleys. The key general integer variables relate to the number of edges used in a block.
The other variables, whilst necessary, flow from the values given to those key variables. This is useful as it enables  the solver we used 
\citep{cplex128} to focus on these key integer variables in its branch and bound search tree. These other variables, being declared as continuous variables, will not be branched on.

\subsection{Heuristic}

It is possible to use our formulation as the basis for a heuristic algorithm for the problem by  only focusing on the first $\tau$ trolleys, instead of focusing on all $T$ trolleys. As discussed above, in our formulation only $z_{ot}$ and $x_{tij}$ need be declared as zero-one, whilst the variables $w_{tb}$ are declared as general integer variables (for all $T$ trolleys).  In the heuristic below we retain the zero-one and integer variables for the first $\tau$ trolleys, but the variables associated with the remaining $T-\tau$ trolleys are
declared as continuous variables. We then solve our formulation to optimality. This will give an assignment of orders to the first $\tau$ trolleys. We then delete those trolleys and the assigned orders from the problem, relabel the orders and trolleys, and repeat. More formally:

\begin{enumerate}[label=(\alph*)]

\item Solve: maximise $\sum_{o \in O}
\sum_{t=1}^\tau b_o z_{ot}$ subject to Equation~(\ref{eq2}) and $ \sum_{o \in O} b_o z_{ot}  \leq  B~t=1,\ldots,T$,
where the $z_{ot}$ are zero-one variables for $t \leq  \tau$, but  continuous variables lying between zero and one for $t >  \tau$. Let the optimal solution be of value $B^*$, then this is the maximum number of baskets that can be allocated to the first $\tau$ trolleys.

\item Solve our formulation with the only integer variables being $z_{ot}$, $x_{tij}$ and  $w_{tb}$,~$t=1,\ldots,\tau$, all other variables being continuous variables. Here we add the constraint $\sum_{o \in O}
\sum_{t=1}^\tau b_o z_{ot} = B^*$
to the formulation so as to make best use of the first $\tau$ trolleys.

\item Delete the orders assigned to trolleys $1,2,\ldots,\tau$ from $O$, delete  the first $\tau$ trolleys and set $T \leftarrow T-\tau$. If $T=0$ or $|O|=0$ go to (e). 

\item Set $\tau \leftarrow \mbox{min}[\tau,T]$. Relabel the remaining trolleys $1,2,\ldots,T$; relabel the orders $1,2,\ldots, |O|$ and go to (a).

\item Optimally route the orders assigned to each trolley using the approach in \cite{valle2017}.

\end{enumerate}
Note here that we relabel the trolleys and orders to ensure that we make use of the constraint that assigns order 1 to trolley 1 (Equation~(\ref{eq37})), as well as the other symmetry constraints (Equations~(\ref{eq38}),(\ref{eq39})).

We would comment in passing here that the general approach adopted in the heuristic algorithm given above is \emph{\textbf{partial integer optimisation}}. By this we mean that from the original formulation a subset of the integer variables are declared as integer, with all of the remaining integer variables being declared as continuous. The resulting mixed-integer problem is then solved,  either to optimality or heuristically. The integer variables are then fixed at the values that they have in this solution and the process repeats with a new subset of integer variables being declared as integer. 

Our partial integer optimisation approach is a \emph{\textbf{matheuristic}}, as it works directly  from a mathematical formulation of the problem under 
consideration~\citep{boschetti2009}.
A related approach to our 
partial integer optimisation approach
is kernel search~\citep{angelelli2012, guastaroba2017}.

As best as we are aware partial integer optimisation  approaches, which essentially lead directly from a  mathematical  formulation to a heuristic, without the necessity of designing a 
problem-specific metaheuristic, have not been significantly explored in the literature.

\section{Improving the linear programming relaxation}
\label{sec:lp}
There are a significant number of constraints (valid inequalities) that can be added to our formulation to improve the value of the linear programming relaxation. In this section we outline these constraints. We also summarise our formulation.

\subsection{Trolley use}

We have $T$ trolleys, each with identical capacity $B$. This means that we know that we must (as a minimum) use at least  $ \lceil \sum_{o \in O} b_o/B \rceil$ trolleys. Since the trolleys are identical we can (arbitrarily) choose to use the first
$ \lceil \sum_{o \in O} b_o/B \rceil$ trolleys.
Hence we can add the constraint:
\begin{optprog}
& \alpha_t & = &  1  & t=1, \ldots, \lceil \sum_{o \in O} b_o/B \rceil  \label{eq40}
\end{optprog}

\subsection{Products picked from a single subaisle}

Suppose that we have a product $p \in P$ such that $|Q(p)|=1$. Then there is just a single (unique) subaisle edge from which this product can be picked. We might reasonably expect that many of the products in $P$ would be of this type, since it is easier, both in terms of picker familiarity with picking, and in terms of product restocking, if a product is stored in just one location/subaisle.

For each subaisle edge $(i,j) \in E$ let $\Phi_{ij}$ be the set of orders such that any order $o \in \Phi_{ij}$ contains at least one product $p 
\in P_o$ with $|Q(p)|=1$ and that product $p$ can only be picked from subaisle edge $(i,j)$.  In other words if we choose to pick any order $o \in \Phi_{ij}$ using trolley $t$ we know that $x_{tij}$ must be one. 
Then the constraint below applies:
\begin{optprog}
& x_{tij} & \geq    (\sum_{o \in \Phi_{ij}} b_o 
z_{ot})/ \mbox{min}[B, \sum_{o \in \Phi_{ij}} b_o]   
& \forall (i,j) \in E, ~t=1,\ldots,T
\end{optprog}
The right-hand side of this constraint involves  summing the number of baskets associated with orders in $\Phi_{ij}$. Clearly if this right-hand side is non-zero then at least one of the orders in $\Phi_{ij}$ is allocated to trolley $t$ and hence $x_{tij}$ should be one. The maximum value that the right-hand side can take is one, as we already have a constraint (Equation~(\ref{eq5})) limiting the number of baskets associated with orders allocated to a trolley to at most $B$.

Since  picking of any order $o \in \Phi_{ij}$ means we must traverse subaisle edge $(i,j)$ we have the constraint:
\begin{optprog}
& \sum_{t=1}^T x_{tij} & \geq &  \lceil \sum_{ o \in \Phi_{ij}} b_o/B \rceil  & \forall (i,j) \in E
\end{optprog}
This constraint ensures that, over all trolleys, each edge $(i,j) \in E$ is traversed enough times to ensure that  it is possible to supply all of the orders which require traversal of that subaisle edge.

\subsection{First and last aisle constraints}
Recall that  $\beta_{1ta}=1$ if the first edge in aisle $a$ in block one, so edge $(a,a+W_A)$, is traversed by trolley $t$, zero otherwise. Then the following constraint applies:
\begin{optprog}
& \beta_{1ta} & \leq & \sum_{e=a}^{W_A} \gamma^{L1}_{te} & a=1, \ldots, W_A, ~ t=1, \ldots, T
\end{optprog}
Here $\beta_{1ta}$ can only be one if some aisle $e$ to the East of aisle $a$, so with $e \geq a$, is the highest indexed aisle visited by trolley $t$  in block one. In a similar fashion we have that:
\begin{optprog}
& \beta_{1ta} & \leq & \sum_{e=1}^{a} \gamma^{F1}_{te} & a=1, \ldots, W_A, ~ t=1, \ldots, T
\end{optprog}
Here $\beta_{1ta}$ can only be one if some aisle $e$ to the West of aisle $a$, so with $e \leq a$, is the lowest indexed aisle visited by trolley $t$  in block one. The corresponding constraints when we have more than one block are:
\begin{optprog}
& \beta_{2ta} & \leq & \sum_{e=a}^{W_A} \gamma^{L2}_{te} & a=1, \ldots, W_A, ~ t=1, \ldots, T \\
& \beta_{2ta} & \leq & \sum_{e=1}^{a} \gamma^{F2}_{te} & a=1, \ldots, W_A, ~ t=1, \ldots, T
\end{optprog}

\subsection{Cross-aisle distance lower bound}
Recall that $ D^{WE}_t$ captures the distance travelled in cross-aisles by trolley $t$. Previously it was defined by Equation~(\ref{eq16}) and Equation~(\ref{eq34}). It is possible to improve the linear programming relaxation estimate of $ D^{WE}_t$ as below by making use of products picked from a single subaisle.

Let $O_1 \subseteq O$ be the set of orders such that each order $o \in O_1$ contains two (or more) products which must be picked from a single subaisle and moreover the aisles in which these subaisles are located are different. 

Now for order $o \in O_1$ identify the lowest indexed aisle in the order such that there is a product in the order which can only be supplied from a subaisle in this aisle (note here that this is not necessarily the lowest indexed aisle in the order, since that aisle might not be associated with such a product). Let this aisle be denoted by $\rho(o)$. Similarly let $\sigma(o)$ be the highest indexed aisle in the order such that there is a product in the order which can only be supplied from a subaisle in this aisle.

Then if order $o$ is allocated to trolley $t$ we know that  aisles $\rho(o)$ and $\sigma(o)$  must be visited by the trolley. This applies even if other orders allocated to the same trolley also require the same products  picked from the single subaisles that were used to identify  $\rho(o)$ and $\sigma(o)$. It is therefore valid to impose the additional constraint:
\begin{optprog}
& D^{WE}_t & \geq    (d_{1\sigma(o)} - d_{1\rho(o)})z_{ot}
& \forall o \in O_1, ~t=1,\ldots,T \label{eq47}
\end{optprog}
This constraint enforces the condition that if order $o \in O_1$ is allocated to trolley $t$ then as a minimum the distance travelled in cross-aisles by trolley $t$ must be $(d_{1\sigma(o)} - d_{1\rho(o)})$. 

Note that this constraint applies even if we have multiple blocks with the products associated with $\rho(o)$ and/or $\sigma(o)$
being picked in blocks $b \geq 2$.

\subsection{Summary}
As we have presented a considerable number of variables and constraints above it is worthwhile to summarise here. Our distance approximation formulation is to optimise~Equation~(\ref{eq36}) subject to 
Equations~(\ref{eq2})-(\ref{eq15}),(\ref{eq21})-(\ref{eq35}),(\ref{eq37})-(\ref{eq47}).
In this formulation the variables  $z_{ot}$ and $x_{tij}$ are declared as zero-one (binary variables) and the variables $w_{tb}$ are declared as
general integer variables,
where $0 \leq w_{tb} \leq \lfloor W_A/2 \rfloor$. The variables $D^{WE}_t$, $z^F_t$ and $z^L_t$ 
 are declared as non-negative continuous variables and the variables $\alpha_{t}$, $\beta_{1ta}$, $\gamma^{F1}_{ta}$, $\gamma^{L1}_{ta}$, $\hat{\alpha_t}$, $\beta_{2ta}$, $\gamma^{F2}_{ta}$, $\gamma^{L2}_{ta}$ and  $y_{tb}$  are declared as continuous variables lying between zero and one.

Solving our formulation yields an allocation of orders to trolleys based on minimising our distance approximation. The picker/trolley routing needed to give a true routing distance can be accomplished either heuristically or optimally. Once orders have been assigned to trolleys then each individual picker/trolley is independent of all others and can be routed separately. In the computational work presented below we used an optimal routing approach based on previous work~\citep{valle2017}.

\section{Computational results}
\label{sec:results}

In this section we present computational results for our distance approximation. We first discuss the test problems that we used and then go on to present computational results for our distance approximation when we have more than one block in the warehouse configuration. We also present results for single block test problems. 

For both the multiple and single block cases we consider 
no-reversal routing, where pickers/trolleys can only reverse direction at the end of a subaisle, as well as reversal routing, where pickers/trolleys can reverse direction after traversing only part of one or more subaisles.

In this section we also discuss the quality of our distance approximation. Results are presented with regard to the effectiveness of our symmetry breaking constraints and the effect of the constraints presented to improve the value of the linear programming relaxation. The results given by our partial integer optimisation heuristic are presented, and compared with a time savings heuristic. Finally results are presented for larger test problems.

We used an Intel Xeon \@ 2.40GHz with 32GB of RAM and Linux as the operating system. The code was written in C++ and Cplex 12.8 \citep{cplex128} was used as the mixed-integer solver.

\subsection{Test problems}
\label{sec:test}

In order to examine the performance of our distance approximation formulation we used the same publicly available test problems as used previously in the branch-and-cut approach of \cite{valle2017}. These test problems  were based on \cite{foodmart} so as to represent a realistic supermarket shopping environment and are  appropriate test problems  given that the  motivation underlying our work is online grocery shopping. A key advantage of using those publicly available test problems  in this paper is that for many of them we know from  \cite{valle2017} the optimal joint batching and routing solution. Hence we are able to make a computational comparison between the results of our approach and those optimal values. 

The test problems considered, as detailed in \cite{valle2017}, are for warehouse configurations involving two blocks (so three cross-aisles). As noted previously above restricting attention to just one block is a situation which is commonly considered in the literature, e.g.~\cite{boysen17, chabot17, gademann01, gademann05, hong2012a, hong2017, koster98, deKoster1999, lu16, menendez17a, menendez17b, menendez17, petersen1997, rao13, zulj18}. Hence in order to examine the performance of our distance approximation approach when we have just a single block we generated additional test problems  for warehouse configurations involving just one block.

The Foodmart database is composed of anonymised customer purchases over two years for a chain of supermarkets. There are a total of 1560 distinct products bought during this period. Orders in Foodmart are generally small in size, with most containing up to 4 or 5 distinct products. As online orders (which inspired the problem considered in this paper) may be composed of dozens of items, we combined different Foodmart orders into a single order. For every customer, all of their purchases made in the first $\Delta$ days are merged into a single order. The merged order may contain not only more distinct products, but also a higher number of units of a single product. A test instance is taken as the $|O|$ orders with the highest number of distinct products. If $|O| = 5$, the 5 largest combined orders make up the instance; if $|O| = 6$, we take the same orders as in $|O| = 5$ plus the sixth largest merged order. We created several test instances for $\Delta = \{5, 10, 20\}$ and $|O| = \{5, \ldots, 15, 20, 25, 30\}$. 

We use a simulated warehouse, 
such as shown in Figure~\ref{fig1}, containing in total $W_A = 8$ aisles. Each aisle contains 33 slots on each side, each slot containing 3 shelves, each shelf holding a sufficient number of one distinct product. This gives a warehouse capacity of 1584 products, enough to hold the 1560 distinct products from the Foodmart database (with some slots left empty). The single block warehouse contains one Northern cross-aisle  and one Southern cross-aisle. The multiple block warehouse contains an extra cross-aisle with 17 slots in the Northern block and 16 slots in the Southern block. Arbitrary values are given, in metres, for aisle lengths and cross-aisle widths, as well as for rack depths and slot widths
(so the pick-face is also of corresponding arbitrary size). The distance from the origin to the vertices in the Northern-most 
cross-aisle is also given. A more complete description of these test problems is given in \cite{valle2017}.

As will be seen later below we also consider larger problems with more orders  ($|O| = \{50, 75\}$) and a warehouse layout containing $W_A = 16$ aisles, with 17 slots on each side, and able to hold a total of 1632 products. Similarly for the multiple block warehouse we added an extra cross-aisle, with 8 slots in the Northern block and 7 slots in the Southern block. For the benefit of future workers all of the test problems that we generated and solved in this paper are publicly available at 
\href{http://www.dcc.ufmg.br/~arbex/orderpicking.html}{http://www.dcc.ufmg.br/$\sim$arbex/orderpicking.html}.

It is worth re-emphasising here the distinct differences between the branch-and-cut approach previously presented \cite{valle2017}
 and the distance approximation approach presented in this paper. These differences are:
\begin{itemize} 
\item
The approach in  \cite{valle2017} is a branch-and-cut approach where cuts (valid inequalities) are added as appropriate during the search tree. The approach presented in this paper is a 
straightforward branch and bound approach where no cuts  are added  at any stage during the search tree. 

\item
The approach in  \cite{valle2017} includes cuts (valid inequalities) associated with eliminating subtours in trolley routes.  The approach presented in this paper does not involve subtours, since it does not directly address the routing problem.

\item 
The approach in  \cite{valle2017} jointly solves the order batching and routing problem, so that it explicitly  decides simultaneously both the batching of orders  and the routes to be adopted so as to optimise over both problems (order batching and routing) directly. To achieve this it makes use of an arc based formulation.

\item  
The approach presented in this paper is an edge based formulation that directly addresses order batching, but uses a distance approximation to influence the batching of orders without directly addressing the routing problem. The advantage of this is that, as will become apparent below, problems  can be solved computationally much faster than the times shown in  \cite{valle2017}.
\end{itemize}

\subsection{Multiple block test problems}

Table~\ref{table2} shows the results for our distance approximation approach when applied to multiple block test problems. In that table the parameter $\Delta$ relates to the labelling of the different test problems in \cite{valle2017} and  $\vert O \vert$ is the total number of orders in each test problem.

In order to provide a direct comparison between the results of our distance approximation approach and an optimal solution approach the final four columns in Table~\ref{table2} show the results from applying the optimal branch-and-cut solution approach in~\cite{valle2017} to these test problems. By contrast with the results presented in~\cite{valle2017} these results are with a more up-to-date version of Cplex (version 12.8) and also include the new symmetry constraints presented above (Equation~(\ref{eq39})).  In those columns T(s) denotes the total computation time in seconds and  a solution value in brackets indicates that the problem terminated at the time limit imposed (6 CPU hours, 21600 seconds) without proving optimality. In such cases the solution value shown corresponds to the value of the best feasible solution known at time limit.

In Table~\ref{table2} we show for each test problem  
the total computation time in seconds to solve our distance approximation formulation to proven optimality. The  optimal distance approximation solution value is shown, as well as the value of the lower bound at the root node of the search tree, so before Cplex commences branching. We also show the total number of tree nodes.

Six columns in Table~\ref{table2} show the results from routing the batching as given by the distance approximation. In those columns we 
show the time taken (in seconds) for routing, the distance associated with the routing solution, and the percentage deviation from the optimal (or best-known) routing solution 
as given in Table~\ref{table2}. 

Recall here that, as discussed above, solving our formulation yields an allocation of orders to trolleys based on minimising our distance approximation. Once orders have been assigned to trolleys then each individual picker/trolley is independent of all others and can be routed separately. The picker/trolley routing needed to give a true routing distance can be accomplished either heuristically or optimally. In the results associated with our distance approximation in Table~\ref{table2} we used an optimal routing approach based on previous work~\citep{valle2017}. So the associated columns in that table show the total time in seconds need to optimally route all pickers/trolleys and the solution value obtained.

To illustrate the results consider the test problem in Table~\ref{table2} with $\Delta=5$ and $|O|=20$ involving $T=3$ trolleys. Our distance approximation formulation was solved to proven optimality in 1.2 seconds, with the optimal solution being 945.9. This required 1838 tree nodes where the initial root node lower bound was 653.9. Given the order batching as decided by our distance approximation then routing of all pickers required 0.2 seconds for routing where no-reversal was enforced, with the associated routing distance being 946.9. 
The associated percentage deviation is zero as the (no-reversal) solution value given in Table~\ref{table2} for this test problem is also 946.9. 
Routing of all pickers required 1.7 seconds for routing where reversal was allowed, with the associated routing distance being 898.4. 
The associated percentage deviation is calculated  using the value of 870.4 taken from the final column in Table~\ref{table2} giving $100(898.4-870.4)/870.4 = 3.22\%$.

\begin{table}[!ht]
\centering

\caption{Multiple blocks results
}
\vspace{0.4cm}

{\scriptsize
\renewcommand{\tabcolsep}{0.75mm} 
\renewcommand{\arraystretch}{1.2} 
\begin{tabular}{|ccc|rrcr|rrr|rrr|rr|rr|}\hline
\multirow{3}{*}{$\Delta$} & \multirow{3}{*}{$|O|$} 
& \multirow{3}{*}{$T$}
& \multicolumn{10}{c|}{Results, distance approximation} & \multicolumn{4}{c|}{Results, optimal branch-and-cut approach} \\
\cline{4-13}
\cline{14-17}
& & & \multicolumn{1}{c}{T(s)} & \multicolumn{1}{c}{Optimal} & \multicolumn{1}{c}{Root node} & \multicolumn{1}{c|}{Number}
 & \multicolumn{3}{c|}{No-reversal routing}  & \multicolumn{3}{c|}{Reversal routing} 
& 
\multicolumn{2}{c|}{No-reversal routing} &  \multicolumn{2}{c|}{Reversal routing} \\ 

& & & &  & \multicolumn{1}{c}{lower} & \multicolumn{1}{c|}{nodes} &   \multicolumn{1}{c}{T(s)} & \multicolumn{1}{r}{Value} & \% dev. & \multicolumn{1}{c}{T(s)} & \multicolumn{1}{r}{Value} & \% dev. 

& \multicolumn{1}{c}{T(s)}& \multicolumn{1}{c|}{Optimal}
& \multicolumn{1}{c}{T(s)}& \multicolumn{1}{c|}{Optimal}\\

& & & &  & \multicolumn{1}{c}{bound} & & & & & & & 
& & \multicolumn{1}{c|}{(best-known)}
& & \multicolumn{1}{c|}{(best-known)}\\

\hline
5	&	5	& 1  &	0.0	&	384.6	&	384.6	&	1	&	0.1	&	384.6	&	0	&	1.0	&	348.6	&	0	 & 0.1 & 384.6 & 1.3 & 348.6 \\
	&	6	& 1 &	0.0	&	384.6	&	384.6	&	1	&	0.1	&	384.6	&	0	&	0.6	&	364.8	&	0	& 0.2 & 384.6 & 0.6 & 364.8 \\
	&	7	& 1 &	0.0	&	384.6	&	384.6	&	1	&	0.1	&	384.6	&	0	&	1.4	&	374.8	&	0	& 0.1 & 384.6 & 1.7 & 374.8 \\
	&	8	& 2 &	0.0	&	543.7	&	453.8	&	17	&	0.3	&	543.7	&	0	&	1.1	&	503.8	&	0	 & 1.3 & 543.7 & 7.2 & 503.8 \\
	&	9	& 2 &	0.0	&	603.2	&	462.3	&	62	&	0.2	&	603.7	&	0.08	&	1.5	&	573.6	&	6.30	 & 1.2 & 603.2 & 8.0 & 539.6 \\
	&	10	& 2 &	0.0	&	611.3	&	493.7	&	31	&	0.1	&	611.4	&	0	&	1.7	&	595.6	&	2.44	 & 1.3 & 611.4 & 8.1 & 581.4 \\
	&	11	& 2 &	0.0	&	641.2	&	540.9	&	44	&	0.2	&	641.7	&	0	&	1.0	&	617.8	&	0.70	& 1.1 & 641.7 & 13.0 & 613.5  \\
	&	12	& 2 &	0.0	&	641.2	&	552.0	&	51	&	0.2	&	641.7	&	0	&	1.1	&	621.8	&	0.06	 & 1.4 & 641.7 & 22.5 & 621.4\\
	&	13	& 2 &	0.0	&	649.3	&	577.0	&	20	&	0.1	&	649.3	&	0	&	1.6	&	641.8	&	2.95	& 1.2 & 649.3 & 14.8 & 623.4  \\
	&	14	& 2 &	0.1	&	691.2	&	592.4	&	85	&	0.2	&	691.7	&	0	&	3.3	&	649.8	&	1.64	& 11.3 & 691.7 & 46.9 & 639.3  \\
	&	15	& 2 &	0.1	&	699.3	&	597.4	&	75	&	0.1	&	699.8	&	0	&	1.1	&	677.3	&	3.66	 & 9.7 & 699.8 & 37.3 & 653.4\\
	&	20	& 3 &	1.2	&	945.9	&	653.9	&	1838	&	0.2	&	946.9	&	0	&	1.7	&	898.4	&	3.22	 & 162.0 & 946.9 & 3035.2 & 870.4\\
	&	25	& 4 &	30.1	&	1154.2	&	721.7	&	22221	&	0.2	&	1155.1	&	0	&	2.1	&	1131.6	&	0.72	 & 13480.1 & 1155.1 & 21600.0 & (1123.5) \\
	&	30	& 4 &	384.5	&	1264.4	&	802.1	&	150557	&	0.2	&	1265.4	&	0	&	1.9	&	1219.9	&	-3.45	 & 21600.0 & (1265.4) & 21600.0 & (1263.5)\\
\multicolumn{3}{|r|}{\textbf{Average}} & \textbf{29.7} & & &  & \textbf{0.2} & 
 
& 	\textbf{0.01}	 
 &
 \textbf{1.5} & 
 
& 	\textbf{1.30}	 & \textbf{2519.4} &  & \textbf{3314.0} & 
\\
\hline
 
10	&	5	& 1 &	0.0	&	384.6	&	384.6	&	1	&	0.2	&	384.6	&	0	&	0.5	&	371.1	&	0	 & 0.2 & 384.6 & 0.4 & 371.1\\
	&	6	& 1 &	0.0	&	384.6	&	384.6	&	1	&	0.2	&	384.6	&	0	&	1.4	&	377.1	&	0	& 0.1 & 384.6 & 1.5 & 377.1 \\
	&	7	& 2 &	0.0	&	613.0	&	561.1	&	20	&	0.2	&	613.0	&	0	&	1.8	&	565.6	&	2.87	 & 2.0 & 613.0 & 6.7 & 549.8 \\
	&	8	& 2 &	0.0	&	681.2	&	553.2	&	98	&	0.3	&	681.2	&	0	&	1.6	&	618.2	&	5.82	 & 3.9 & 681.2 & 7.1 & 584.2 \\
	&	9	& 2 &	0.0	&	681.2	&	589.4	&	22	&	0.1	&	681.2	&	0	&	3.9	&	637.4	&	0	 & 2.1 & 681.2 & 55.0 & 637.4\\
	&	10	& 2 &	0.0	&	729.2	&	640.8	&	49	&	0.1	&	729.7	&	0.05	&	3.6	&	695.7	&	5.12	 & 6.8 & 729.3 & 63.9 & 661.8\\
	&	11	& 2 &	0.0	&	729.2	&	687.2	&	13	&	0.1	&	729.7	&	0	&	2.8	&	719.7	&	2.84	 & 3.4 & 729.7 & 655.0 & 699.8\\
	&	12	& 2 &	0.0	&	731.2	&	697.5	&	7	&	0.2	&	731.3	&	0	&	2.8	&	723.8	&	2.27	 & 2.6 & 731.3 & 39.5 & 707.7\\
	&	13	& 2 &	0.1	&	769.2	&	644.9	&	74	&	0.4	&	769.2	&	0	&	2.7	&	727.7	&	0.28	& 74.8 & 769.2 & 497.9 & 725.7  \\
	&	14	& 2 &	0.1	&	769.2	&	673.2	&	42	&	0.4	&	769.2	&	0	&	6.7	&	749.3	&	2.95	 & 47.4 & 769.2 & 389.3 & 727.8\\
	&	15	& 3 &	0.3	&	929.8	&	704.0	&	558	&	0.4	&	930.0	&	0	&	5.3	&	908.5	&	2.93	& 59.7 & 930.0 & 664.6 & 882.6 \\
	&	20	& 3 &	0.5	&	1027.8	&	767.6	&	807	&	0.2	&	1028.3	&	0	&	3.4	&	1007.3	&	1.50	 & 153.0 & 1028.3 & 10923.4 & 992.4 \\
	&	25	& 4 &	11.5	&	1264.4	&	868.2	&	7276	&	0.9	&	1264.6	&	0.01	&	2.7	&	1233.7	&	-2.56	 & 11768.2 & 1264.5 & 21600.0 & (1266.1)\\
	&	30	& 4 &	137.5	&	1362.4	&	975.5	&	39598	&	0.4	&	1362.9	&	0.04	&	4.9	&	1314.1	&	-2.34	& 21600.0 & (1362.4) & 21600.0 & (1345.6) \\
\multicolumn{3}{|r|}{\textbf{Average}} & \textbf{10.7} & & &  & \textbf{0.3} &  
 
& 	\textbf{0.01}	 

& \textbf{3.2} & 
 
& 	\textbf{1.55}	 
 & \textbf{2408.9} &  & \textbf{4036.0} & 
\\
\hline
 
20	&	5	& 2 &	0.0	&	623.3	&	572.6	&	10	&	0.1	&	623.3	&	0	&	2.1	&	573.8	&	0	& 1.3 & 623.3 & 7.7 & 573.8 \\
	&	6	& 2 &	0.0	&	729.2	&	697.2	&	7	&	0.2	&	729.3	&	0	&	1.1	&	668.2	&	1.83	 & 2.0 & 729.3 & 20.6 & 656.2\\
	&	7	& 2 &	0.0	&	769.2	&	729.7	&	9	&	0.3	&	769.2	&	0	&	1.5	&	719.7	&	4.33	 & 8.7 & 769.2 & 14.0 & 689.8 \\
	&	8	& 2 &	0.0	&	769.2	&	740.7	&	3	&	0.1	&	769.2	&	0	&	2.3	&	737.7	&	5.72	& 12.0 & 769.2 & 26.1 & 697.8  \\
	&	9	& 2 &	0.0	&	769.2	&	709.2	&	19	&	0.3	&	769.2	&	0	&	3.5	&	741.7	&	1.92	& 12.0 & 769.2 & 115.6 & 727.7 \\
	&	10	& 3 &	0.2	&	989.8	&	784.0	&	452	&	0.4	&	989.9	&	0	&	2.7	&	946.5	&	2.82	 & 19.7 & 989.9 & 179.1 & 920.5\\
	&	11	& 3 &	0.3	&	1075.8	&	804.7	&	1098	&	0.3	&	1075.8	&	0	&	3.4	&	1004.5	&	2.45	& 82.9 & 1075.8 & 309.2 & 980.5 \\
	&	12	& 3 &	0.4	&	1113.8	&	797.0	&	1295	&	0.2	&	1114.0	&	0.02	&	3.8	&	1034.4	&	3.00	& 506.5 & 1113.8 & 871.0 & 1004.3 \\
	&	13	& 3 &	0.4	&	1113.8	&	848.5	&	1298	&	0.6	&	1114.0	&	0.02	&	3.5	&	1025.0	&	1.58	 & 676.7 & 1113.8 & 852.7 & 1009.1\\
	&	14	& 3 &	1.0	&	1113.8	&	915.8	&	711	&	0.5	&	1113.9	&	0	&	3.1	&	1068.2	&	5.65	& 421.4 & 1113.9 & 506.1 & 1011.1 \\
	&	15	& 3 &	0.9	&	1115.8	&	936.2	&	622	&	0.8	&	1115.9	&	0	&	4.1	&	1068.5	&	3.87	 & 348.9 & 1115.9 & 2308.9 & 1028.7 \\
	&	20	& 4 &	24.8	&	1438.2	&	1184.2	&	7722	&	0.3	&	1438.8	&	0.01	&	3.4	&	1387.3	&	1.00	& 5752.8 & 1438.7 & 21600.0 & (1373.5) \\
	&	25	& 5 &	221.8	&	1727.0	&	1289.6	&	118215	&	0.3	&	1727.2	&	-0.12	&	3.4	&	1672.4	&	-1.18	& 21600.0 & (1729.2) & 21600.0 & (1692.3) \\
	&	30	& 6 &	10639.7	&	1995.4	&	1344.0	&	3808495	&	0.5	&	1996.2	&	-3.20	&	4.8	&	1919.1	&	-1.32	 & 21600.0 & (2062.2) & 21600.0 & (1944.7)\\
\multicolumn{3}{|r|}{\textbf{Average}} & \textbf{777.8} & & &  & \textbf{0.4} &  
 
& 	\textbf{-0.23}	 
& \textbf{3.1} & 
 
& 	\textbf{2.26}	 & \textbf{3646.1} &  & \textbf{5000.8} & 
\\
\hline
\multicolumn{3}{r}{\textbf{Average all}} & \textbf{272.8} & & &  \multicolumn{1}{r}{} & \multicolumn{1}{r}{\textbf{0.3}} &  

& 	 \multicolumn{1}{r}{\textbf{-0.07}}	 
& \multicolumn{1}{r}{\textbf{2.6}}

& 	& \multicolumn{1}{r}{\textbf{1.70}}	 &
\multicolumn{1}{r}{\textbf{2858.1}} & \multicolumn{1}{r}{} & \multicolumn{1}{r}{\textbf{4116.9}} & \multicolumn{1}{r}{}
\\
\end{tabular}
}

\label{table2}
\end{table}

Considering Table~\ref{table2} the key issues for our distance approximation approach are the quality of the final routing obtained and the computation time required, as compared with the values seen in the last four columns of Table~\ref{table2} for the optimal  approach.

For no-reversal routing 
(where pickers/trolleys can only reverse direction at the end of a subaisle)
we can see that for 33 of the 42  test problems we obtain the same solution as the optimal approach (as indicated by a percentage deviation of zero). For two of the larger problems that are not solved to proven optimality  we obtain better solutions (as indicated by negative percentage deviations). Over all the test problems considered the average percentage deviation is 
$-0.07$\%, so a very slight improvement in routing quality. However the computation time required is much less. For no-reversal routing the average computation time 
in Table~\ref{table2} 
 is 272.8 seconds to solve the distance approximation, and a further 0.3 seconds for routing, so a total of 273.1 seconds. By comparison the average computation time for the optimal  approach for no-reversal routing  is 
2858.1 seconds. So for the test problems examined the distance approximation approach presented in this paper requires only $100(273.1/2858.1)= 9.6\%$ of the time required by the optimal  approach, but gives very slightly better  quality solutions when considering no-reversal routing.

With respect to reversal routing
(where pickers/trolleys can reverse direction after traversing only part of one or more subaisles) it is clear from the distance approximation developed above that at no stage did we attempt to account for reversal in a subaisle. For this reason we would expect the quality of the routing to be poorer than for the no-reversal case. However there may be computational advantages in using our distance approximation approach to batch orders for reversal routing.

Considering Table~\ref{table2} we can see that the percentage deviations values are worse than for no-reversal routing. 
For five of the larger problems that are not solved to proven optimality by the optimal approach we obtain better solutions (as indicated by negative percentage deviations). 
Over all the test problems considered the average percentage deviation for reversal routing was 1.70\%. The average computation time was 272.8+2.6=275.4 seconds. This compares with an average computation time of 4116.9 seconds for the optimal approach.  So for the test problems examined the distance approximation approach presented in this paper requires only $100(275.2/4116.9)= 6.7\%$ of the time required by the optimal approach, but gives solutions on average 1.70\% worse when considering reversal routing. Obviously there is a value judgement to be made here as to whether the significantly lower computation time outweighs the slightly longer routing distance.

\subsection{Single  block test problems}

Table~\ref{table3} gives the results for the single block test problems. This table has the same format as Table~\ref{table2}. As for Table~\ref{table2} the key issues for our distance approximation approach are the quality of the final routing obtained and the computation time required, as compared with the values seen in  the last four columns of Table~\ref{table3} for the optimal  approach.

We stated above that for single block instances with no-reversal routing our distance approximation approach gives an optimal allocation of orders to trolleys and an optimal routing.
This can be seen in  Table~\ref{table3} where, 
 for all  of the test problems, the optimal distance approximation solution value is equal to the no-reversal routing value and we obtain the same solution as the optimal branch-and-cut approach (as indicated by a percentage deviation of zero).
 However the computation time required by our distance approximation approach is significantly less than that for the optimal 
branch-and-cut approach.
For 
no-reversal routing the average computation time in 
Table~\ref{table3} is 0.4+0.2=0.6 seconds.
By comparison the average computation time for no-reversal routing for the optimal approach is 
2807.0 seconds, so a factor of approximately 4700 times more. 

Considering Table~\ref{table3} then over all the test problems considered the average percentage deviation for reversal routing was 2.66\%. The average computation time was 0.4+1.5=1.9 seconds.
By comparison the average computation time for reversal routing 
for the optimal approach  is 
3740.9 seconds, so a factor of approximately 2000 times more.
So here, as before for reversal routing with multiple blocks, there is a value judgement to be made as to whether the much lower computation time outweighs the slightly longer routing distance.

In the description of reversal and no-reversal routing given above in association with Figure~\ref{fig2} and Figure~\ref{fig3} we noted that allowing routes involving reversal may involve less distance than no-reversal routing. As can be seen from the results in Table~\ref{table3}, both  for our distance approximation and for the optimal branch-and-cut approach, there may be no gain in allowing reversal routing (e.g.~consider the second problem in Table~\ref{table3} with $\Delta=5$, $|O|=6$ where the routing distance is 360.6). For the single block instances and 
the distance approximation approach
in Table~\ref{table3} there is no gain from reversal routing in 14 of the 42 cases. By contrast for the multiple block instances and 
the distance approximation approach
in Table~\ref{table2} there is a gain from reversal routing in all 42 cases.

\begin{table}[!ht]
\centering

\caption{Single block results
}
\vspace{0.4cm}

{\scriptsize
\renewcommand{\tabcolsep}{0.75mm} 
\renewcommand{\arraystretch}{1.2} 
\begin{tabular}{|ccc|rrcr|rrr|rrr|rr|rr|}\hline
\multirow{3}{*}{$\Delta$} & \multirow{3}{*}{$|O|$} 
& \multirow{3}{*}{$T$}
& \multicolumn{10}{c|}{Results, distance approximation} & \multicolumn{4}{c|}{Results, optimal branch-and-cut  approach} \\
\cline{4-13}
\cline{14-17}
& & & \multicolumn{1}{c}{T(s)} & \multicolumn{1}{c}{Optimal} & \multicolumn{1}{c}{Root node} & \multicolumn{1}{c|}{Number}
 & \multicolumn{3}{c|}{No-reversal routing}  & \multicolumn{3}{c|}{Reversal routing} 
& 
\multicolumn{2}{c|}{No-reversal routing} &  \multicolumn{2}{c|}{Reversal routing} \\ 

& & & &  & \multicolumn{1}{c}{lower} & \multicolumn{1}{c|}{nodes} &   \multicolumn{1}{c}{T(s)} & \multicolumn{1}{r}{Value} & \% dev. & \multicolumn{1}{c}{T(s)} & \multicolumn{1}{r}{Value} & \% dev. 

& \multicolumn{1}{c}{T(s)}& \multicolumn{1}{c|}{Optimal}
& \multicolumn{1}{c}{T(s)}& \multicolumn{1}{c|}{Optimal}\\

& & & &  & \multicolumn{1}{c}{bound} & & & & & & & 
& & \multicolumn{1}{c|}{(best-known)}
& & \multicolumn{1}{c|}{(best-known)}\\

\hline
 
5	&	5	& 1  &	0.0	&	360.6	&	360.6	&	1	&	0.1	&	360.6	&	0	&	0.2	&	338.6	&	0	 & 0.0 & 360.6 & 0.2 & 338.6 \\
	&	6	& 1  &	0.0	&	360.6	&	360.6	&	1	&	0.1	&	360.6	&	0	&	0.5	&	360.6	&	0	 & 0.1 & 360.6 & 0.3 & 360.6\\
	&	7	& 1  &	0.0	&	360.6	&	360.6	&	1	&	0.0	&	360.6	&	0	&	0.6	&	360.6	&	0	 & 0.1 & 360.6 & 0.5 & 360.6  \\
	&	8	& 2  &	0.0	&	547.2	&	547.2	&	1	&	0.1	&	547.2	&	0	&	0.5	&	531.7	&	0.45	 & 0.3 & 547.2 & 4.5 & 529.3  \\
	&	9	& 2 &	0.0	&	629.2	&	555.3	&	9	&	0.1	&	629.2	&	0	&	0.7	&	609.8	&	8.56	 & 0.4 & 629.2 & 4.5 & 561.7  \\
	&	10	& 2 &	0.0	&	639.2	&	639.2	&	1	&	0.1	&	639.2	&	0	&	0.7	&	597.7	&	0	 & 0.5 & 639.2 & 4.3 & 597.7 \\
	&	11	& 2 &	0.0	&	711.2	&	711.2	&	1	&	0.1	&	711.2	&	0	&	0.5	&	673.2	&	5.90	& 1.5 & 711.2 & 5.9 & 635.7   \\
	&	12	& 2 &	0.0	&	711.2	&	671.8	&	3	&	0.1	&	711.2	&	0	&	0.5	&	673.2	&	2.36	& 1.8 & 711.2 & 9.4 & 657.7  \\
	&	13	& 2 &	0.0	&	711.2	&	711.2	&	1	&	0.1	&	711.2	&	0	&	0.7	&	673.2	&	0.81	 & 1.3 & 711.2 & 8.1 & 667.8    \\
	&	14	& 2 &	0.0	&	719.3	&	683.3	&	7	&	0.1	&	719.3	&	0	&	1.3	&	685.8	&	2.10	 & 4.9 & 719.3 & 23.8 & 671.7  \\
	&	15	& 2 &	0.0	&	721.2	&	721.2	&	1	&	0.2	&	721.2	&	0	&	0.9	&	713.2	&	3.71	 & 3.8 & 721.2 & 101.6 & 687.7 \\
	&	20	& 3 &	0.1	&	987.9	&	768.4	&	118	&	0.2	&	987.9	&	0	&	1.9	&	956.4	&	3.08	 & 38.0 & 987.9 & 878.6 & 927.8   \\
	&	25	& 4 &	1.7	&	1268.4	&	870.4	&	1922	&	0.3	&	1268.4	&	0	&	1.2	&	1252.9	&	6.10	  & 21600.0 & (1268.4) & 21600.0 & (1180.9)   \\
	&	30	& 4 &	6.6	&	1358.5	&	1084.4	&	2588	&	0.4	&	1358.5	&	0	&	1.5	&	1297.1	&	2.21	 & 21600.0 & (1358.5) & 21600.0 & (1269.1)  \\
\multicolumn{3}{|r|}{\textbf{Average}} & \textbf{0.6} & & &  & \textbf{0.1} &  
 
& 	\textbf{0}	 
& \textbf{0.8} & 
 
& 	\textbf{2.52}	 
& \textbf{3089.5} &  & \textbf{3160.1} &    
\\
\hline
 
10	&	5	& 1 &	0.0	&	360.6	&	360.6	&	1	&	0.1	&	360.6	&	0	&	0.8	&	360.6	&	0	 & 0.1 & 360.6 & 0.2 & 360.6     \\
	&	6	& 1 &	0.0	&	360.6	&	360.6	&	1	&	0.0	&	360.6	&	0	&	0.6	&	360.6	&	0	 & 0.0 & 360.6 & 0.3 & 360.6  \\
	&	7	& 2 &	0.0	&	639.2	&	566.2	&	7	&	0.1	&	639.2	&	0	&	2.1	&	633.2	&	11.68	  & 0.5 & 639.2 & 5.1 & 567.0    \\
	&	8	& 2 &	0.0	&	639.2	&	639.2	&	1	&	0.1	&	639.2	&	0	&	2.1	&	639.2	&	5.53	 & 0.6 & 639.2 & 4.7 & 605.7   \\
	&	9	& 2 &	0.0	&	639.2	&	639.2	&	1	&	0.1	&	639.2	&	0	&	2.9	&	639.2	&	0	& 0.4 & 639.2 & 15.5 & 639.2  \\
	&	10	& 2 &	0.0	&	721.2	&	721.2	&	1	&	0.2	&	721.2	&	0	&	1.2	&	719.3	&	5.90	 & 2.1 & 721.2 & 35.5 & 679.2 \\
	&	11	& 2 &	0.0	&	721.2	&	721.2	&	1	&	0.1	&	721.2	&	0	&	1.4	&	719.3	&	5.59	& 1.3 & 721.2 & 8.6 & 681.2     \\
	&	12	& 2 &	0.0	&	721.2	&	721.2	&	1	&	0.1	&	721.2	&	0	&	1.8	&	721.2	&	2.85	& 1.2 & 721.2 & 15.1 & 701.2   \\
	&	13	& 2 &	0.0	&	721.2	&	721.2	&	1	&	0.1	&	721.2	&	0	&	1.7	&	721.2	&	1.12	 & 4.5 & 721.2 & 734.8 & 713.2   \\
	&	14	& 2 &	0.0	&	721.2	&	721.2	&	1	&	0.2	&	721.2	&	0	&	1.9	&	721.2	&	0.84	 & 3.9 & 721.2 & 566.3 & 715.2 \\
	&	15	& 3 &	0.1	&	979.8	&	754.7	&	114	&	0.2	&	979.8	&	0	&	2.1	&	917.8	&	1.53	 & 24.4 & 979.8 & 345.3 & 904.0  \\
	&	20	& 3 &	0.0	&	999.8	&	966.4	&	6	&	0.2	&	999.8	&	0	&	2.9	&	999.8	&	1.63	 & 8.5 & 999.8 & 1354.5 & 983.8 \\
	&	25	& 4 &	0.6	&	1278.4	&	1038.5	&	535	&	0.1	&	1278.4	&	0	&	2.0	&	1278.4	&	6.49	 & 7536.1 & 1278.4 & 21600.0 & (1200.5)\\
	&	30	& 4 &	1.1	&	1350.4	&	1176.5	&	240	&	0.5	&	1350.4	&	0	&	4.3	&	1318.4	&	3.12	 & 21600.0 & (1350.4) & 21600.0 & (1278.5)   \\
\multicolumn{3}{|r|}{\textbf{Average}} & \textbf{0.1} & & &  & \textbf{0.2} &  
 
& 	\textbf{0}	 
& \textbf{2.0} & 
 
& 	\textbf{3.31}	 & \textbf{2084.5} &  & \textbf{3306.1} &   
\\
\hline
 
20	&	5	& 2 &	0.0	&	711.2	&	711.2	&	1	&	0.1	&	711.2	&	0	&	0.7	&	631.7	&	1.33	 & 0.3 & 711.2 & 6.2 & 623.4  \\
	&	6	& 2 &	0.0	&	721.2	&	721.2	&	1	&	0.2	&	721.2	&	0	&	0.6	&	703.7	&	3.00	 & 0.4 & 721.2 & 4.9 & 683.2\\
	&	7	& 2 &	0.0	&	721.2	&	721.2	&	1	&	0.2	&	721.2	&	0	&	0.9	&	711.2	&	2.30	& 0.5 & 721.2 & 8.0 & 695.2   \\
	&	8	& 2 &	0.0	&	721.2	&	721.2	&	1	&	0.2	&	721.2	&	0	&	0.8	&	721.2	&	3.74	& 0.6 & 721.2 & 8.9 & 695.2      \\
	&	9	& 2 &	0.0	&	721.2	&	721.2	&	1	&	0.1	&	721.2	&	0	&	1.1	&	721.2	&	1.68	& 2.2 & 721.2 & 22.9 & 709.3  \\
	&	10	& 3 &	0.0	&	999.8	&	998	&	6	&	0.1	&	999.8	&	0	&	1.3	&	953.9	&	1.26	 & 2.7 & 999.8 & 44.7 & 942.0  \\
	&	11	& 3 &	0.0	&	1071.8	&	974.7	&	14	&	0.3	&	1071.8	&	0	&	0.7	&	1031.8	&	3.10	  & 13.5 & 1071.8 & 249.8 & 1000.8  \\
	&	12	& 3 &	0.0	&	1071.8	&	979.5	&	16	&	0.2	&	1071.8	&	0	&	1.1	&	1027.8	&	1.97	 & 23.5 & 1071.8 & 160.1 & 1007.9 \\
	&	13	& 3 &	0.0	&	1071.8	&	1035.8	&	10	&	0.4	&	1071.8	&	0	&	1.8	&	1040.3	&	1.80	  & 16.1 & 1071.8 & 262.7 & 1021.9 \\
	&	14	& 3 &	0.0	&	1071.8	&	1071.8	&	1	&	0.2	&	1071.8	&	0	&	1.9	&	1045.8	&	1.70	& 20.8 & 1071.8 & 477.1 & 1028.3   \\
	&	15	& 3 &	0.0	&	1081.8	&	1081.8	&	1	&	0.2	&	1081.8	&	0	&	1.7	&	1081.8	&	4.19	 & 57.5 & 1081.8 & 544.6 & 1038.3   \\
	&	20	& 4 &	0.1	&	1432.4	&	1396.4	&	3	&	0.2	&	1432.4	&	0	&	3.1	&	1412.9	&	0.88	 & 2117.8 & 1432.4 & 21600.0 & (1400.6)   \\
	&	25	& 5 &	0.3	&	1791.0	&	1647.0	&	198	&	0.4	&	1791.0	&	0	&	2.7	&	1729.6	&	1.91	 & 21600.0 & (1791.0) & 21600.0 & (1697.2)  \\
	&	30	& 6 &	8.0	&	2069.7	&	1759.5	&	4504	&	0.3	&	2069.7	&	0	&	3.0	&	2016.2	&	1.29	 & 21600.0 & (2069.7) & 21600.0 & (1990.6)  \\
\multicolumn{3}{|r|}{\textbf{Average}} & \textbf{0.6} & & &  & \textbf{0.2} &  
 
& 	\textbf{0}	 
& \textbf{1.5} & 
 
& 	\textbf{2.15}	 & \textbf{3246.9} &  & \textbf{4756.4} &    \\
\hline
\multicolumn{3}{r}{\textbf{Average all}} & \textbf{0.4} & & &  \multicolumn{1}{r}{} & \multicolumn{1}{r}{\textbf{0.2}} &  

& 	 \multicolumn{1}{r}{\textbf{0}}	 
& \multicolumn{1}{r}{\textbf{1.5}}

& 	& \multicolumn{1}{r}{\textbf{2.66}}	 
& \multicolumn{1}{r}{\textbf{2807.0}} & \multicolumn{1}{r}{} & \multicolumn{1}{r}{\textbf{3740.9}} & \multicolumn{1}{r}{} 
\\
\end{tabular}

}
\label{table3}
\end{table}

With regard to our single block test problems then, as   noted previously above, restricting attention to just one block is a situation which is commonly considered in the literature, 
e.g.~\cite{boysen17, chabot17, gademann01, gademann05, hong2012a, hong2017, koster98, deKoster1999, 
lu16, menendez17a, menendez17b, menendez17, petersen1997, rao13, zulj18}. Whilst Table~\ref{table3} has presented a comparison between our distance approximation approach and an optimal branch-and-cut approach there is a question as to how well our distance approximation approach compares with these other published approaches to single block problems.

In our view it is clearly impracticable to attempt our own implementation of these other approaches.
 Instead  we implemented a time savings based heuristic to provide a basis for evaluating the performance of our distance approximation approach as compared with a standard  heuristic approach drawn from the literature. In support of this strategy of comparing our results with the results from a time savings heuristic we would stress here that making use of a time savings based heuristic in order to provide a basis for evaluating the performance of another heuristic is common in the literature, e.g.~see~\cite{albareda2009, deKoster1999, 
henn2012_2, henn2012, hong2017, matusiak2017, gils18a,  zulj18}.

For space reasons we do not present here detailed results,  rather Table~\ref{table_single_extra} shows average computation times and average percentage deviations for the single block test problems in Table~\ref{table3}.  The computation time for batching  for our distance approximation  approach does not depend on the routing approach (no-reversal or reversal) adopted. However the computation time for the time savings heuristic does depend upon the routing approach adopted.
For the time savings heuristic for the single  block test problems 
with $\Delta=20$ in Table~\ref{table_single_extra} 
the average computation time for no-reversal routing was 6.9 seconds, with the average percentage deviation being 6.38\%, whilst for reversal routing the average computation time was 17.3 seconds with the average percentage deviation being 3.15\%. 

Examining 
Table~\ref{table_single_extra}
we can see that for the single block test problems  our distance approximation   significantly outperforms the time savings heuristic, both with regard to computation time and routing quality. For these instances the average computation times for our distance approximation are very small. 

\begin{table}[!ht]
\centering

\caption{Single block: comparison with a time savings based heuristic}
\vspace{0.4cm}

{\scriptsize
\renewcommand{\tabcolsep}{2.9mm} 
\renewcommand{\arraystretch}{1.1} 
\begin{tabular}{|ccrcc|}
\hline
 $\Delta$ & Case & T(s) & No-reversal  & Reversal  \\
 & & & routing \% dev. & routing \% dev. \\
\hline
5 & Distance approximation, Table~\ref{table3} & 0.6& 0 & 2.52 \\
& Time savings & 7.1, 7.1 & 3.94 & 5.54 \\

10 & Distance approximation, Table~\ref{table3} & 0.1 & 0 & 3.31 \\

& Time savings & 7.3, 11.9 & 3.87 & 13.57 \\

20 & Distance approximation, Table~\ref{table3} & 0.6 & 0 & 2.15\\
& Time savings & 6.9, 17.3 & 6.38 & 3.15 \\     

\hline
\end{tabular}
}
\label{table_single_extra}
\end{table}

\subsection{Distance approximation quality}

Above we have compared the results obtained from our distance approximation approach with  the optimal  branch-and-cut approach presented in \cite{valle2017}.
 However a further issue of interest is the quality of our distance approximation. In other words how close are the values of the optimal distance approximation solution and the routing solution based on the order batching given by the distance approximation. To examine this issue we computed, from Table~\ref{table2} and Table~\ref{table3}, the average value of 100(Routing solution value - Optimal distance approximation value)/(Optimal distance approximation value). 

For  Table~\ref{table2} the average value for this measure was 0.03\% for no-reversal routing and $-4.49$\% for reversal routing. 
For  Table~\ref{table3} the average value for this measure was 0\% for no-reversal routing and $-2.37$\% for reversal routing.

Hence we can conclude here that for no-reversal routing (for the test problems examined) our distance approximation very closely approximates the underlying routing distance. This is important as it indicates that any algorithm aiming to minimise our distance approximation will, by default, also minimise the underlying routing distance that results from the batching decided.

For reversal routing the negative values obtained indicate that (as we might reasonably expect) our distance approximation over-estimates the routing distance when reversal is allowed. However, as the computation time comparison above between 
the optimal approach and our distance approximation approach 
 indicated, there may well be computational advantages in adopting a distance approximation to guide order batching so as to lessen the computational effort required as compared to trying to achieve  optimal joint batching and routing, especially as problem size increases.

\subsection{Symmetry}

To illustrate the importance of imposing our new symmetry breaking constraint (Equation~(\ref{eq39})) we solved all of the test problems considered in Table~\ref{table2} and
Table~\ref{table3}, but without imposing that constraint. For space reasons however we only report summary results here.

Table~\ref{table4} shows the results for the branch-and-cut approach of \cite{valle2017} both with Equation~(\ref{eq39}) 
(as in Table~\ref{table2} and 
Table~\ref{table3}) and without Equation~(\ref{eq39}). In that table we show for each value of $\Delta$ the  number of problems unsolved at time limit (6 CPU hours), as well as the average time (in seconds). In general we can say that  imposing Equation~(\ref{eq39}) improves the results, in particular for no-reversal routing.

Table~\ref{table5} shows the same information as Table~\ref{table4}, but for the distance approximation approach presented in this paper. The symmetry breaking results in Table~\ref{table5} are taken from Table~\ref{table2} and Table~\ref{table3}. We can see that  including Equation~(\ref{eq39}) has a significant effect in terms of reducing the computation time, especially here for the multiple blocks case.

\begin{table}[!ht]

\caption{Results with and without symmetry breaking constraint Equation~(\ref{eq39}), \cite{valle2017}}
\vspace{0.4cm}

\centering
{\scriptsize
\renewcommand{\tabcolsep}{0.75mm} 
\renewcommand{\arraystretch}{1.2} 
\begin{tabular}{|cl|rr|rr|rr|rr|}\hline
\multirow{4}{*}{$\Delta$} 	&		&	\multicolumn{4}{c|}{Multiple blocks}	&		\multicolumn{4}{c|}{Single block} \\	\cline{3-10}	   
	&		&		\multicolumn{2}{c|}{No-reversal}	&	\multicolumn{2}{c|}{Reversal}	&	\multicolumn{2}{c|}{No-reversal}	&	\multicolumn{2}{c|}{Reversal}		\\	   
	&		&	Symmetry 	&	No symmetry 	&	Symmetry 	&	No symmetry 	&	Symmetry 	&	No symmetry 	&	Symmetry 	&	No symmetry 	\\	   
	&		&	 breaking	&	 breaking	&	 breaking	&	breaking	&	 breaking	&	 breaking	&	 breaking	&	 breaking		\\

5	&	Number unsolved	&	1	&	2	&	2	&	2	&	2	&	2	&	2	&	2		\\	   
	&	Average time (s)	&	2519.4	&	3119.5	&	3314.0	&	3607.0	&	3089.5	&	3089.7	&	3160.1	&	3241.3		\\	   
10	&	Number unsolved	&	1	&	2	&	2	&	2	&	1	&	2	&	2	&	2		\\	   
	&	Average time (s)	&	2408.9	&	3124.2	&	4036.0	&	4452.4	&	2084.5	&	3091.2	&	3306.1	&	3419.9		\\	   
20	&	Number unsolved	&	2	&	3	&	3	&	3	&	2	&	2	&	3	&	3		\\	   
	&	Average time (s)	&	3646.1	&	4778.2	&	5000.8	&	5119.0	&	3246.9	&	4314.1	&	4756.4	&	4893.4		\\	   
Average	&	Number unsolved	&	1.3	&	2.3	&	2.3	&	2.3	&	1.7	&	2.0	&	2.3	&	2.3		\\	   
	&	Average time (s)	&	2858.1	&	3674.0	&	4116.9	&	4392.8	&	2807.0	&	3498.3	&	3740.9	&	3851.5		\\	 
 
\hline
\end{tabular}
}

\label{table4}
\end{table}

\begin{table}[!ht]
\centering

\caption{Results with and without symmetry breaking constraint Equation~(\ref{eq39}), distance approximation}

\vspace{0.4cm}

{\scriptsize
\renewcommand{\tabcolsep}{0.75mm} 
\renewcommand{\arraystretch}{1.2} 
\begin{tabular}{|cl|rr|rr|}\hline
\multirow{4}{*}{$\Delta$} 
	&		&	\multicolumn{2}{c|}{Multiple blocks}		&		\multicolumn{2}{c|}{Single block}			\\	   
	&		&	Symmetry 	&	No symmetry 	&	Symmetry 	&	No symmetry 	\\	
	&		&	breaking	&	breaking	&	 breaking	&	breaking	\\	      
5	&	Number unsolved	&	0	&	0	&	0	&	0		\\	   
	&	Average time (s)	&	29.7	&	75.9	&	0.6	&	1.4		\\	   
10	&	Number unsolved	&	0	&	0	&	0	&	0		\\	   
	&	Average time (s)	&	10.7	&	31.3	&	0.1	&	0.2		\\	   
20	&	Number unsolved	&	0	&	1	&	0	&	0		\\	   
	&	Average time (s)	&	777.8	&	1589.8	&	0.6	&	1.2		\\	   
Average	&	Number unsolved	&	0	&	0.3	&	0	&	0		\\	   
	&	Average time (s)	&	272.8	&	565.7	&	0.4	&	0.9		\\	 
\hline
\end{tabular}
}

\label{table5}
\end{table}

\subsection{Linear programming relaxation improvement}

To illustrate the effect of the constraints, Equations~(\ref{eq40})-(\ref{eq47}), presented above which improve  the value of the linear programming relaxation of our formulation Table~\ref{table_lp} 
presents, for the problems considered previously above,
the percentage improvement in the root node lower bound (RNLB), as measured by 100((RNLB with  Equations~(\ref{eq40})-(\ref{eq47})) - (RNLB without  Equations~(\ref{eq40})-(\ref{eq47})))/(RNLB without  Equations~(\ref{eq40})-(\ref{eq47})). 

For the smaller problems our additional constraints make little difference. However for the larger problems substantial improvements in the RNLB were observed, both for multiple block and single block problems. This behaviour can be seen in Table~\ref{table_lp}. For example for the case $\Delta=20$ the average percentage improvement in the RNLB for the multiple blocks problems was 19.36\% across all values of $|O|$, but 43.39\% for the four problems with 	$|O|=15,20,25,30$.

\begin{table}[!ht]
\centering

\caption{Root node, average lower bound percentage improvement}
\vspace{0.4cm}

{\scriptsize
\renewcommand{\tabcolsep}{2.9mm} 
\renewcommand{\arraystretch}{1.2} 
\begin{tabular}{|c|cc|cc| }
\hline
$\Delta$ & \multicolumn{2}{c|}{Single block} & \multicolumn{2}{c|}{Multiple blocks}\\
  & All $|O|$ & $|O|=15,20,25,30$ & All $|O|$ & $|O|=15,20,25,30$ \\
 \hline
5	&	12.97	&	34.11	&	5.76	&		16.21	\\
10	&	14.10	&	37.04		&	9.36	&	22.90	\\
20	&	12.38	&	31.48		&	19.36	&	43.39	\\
\hline
\multicolumn{1}{r}{\textbf{Average all}}   & \multicolumn{1}{c}{\textbf{    13.15}} & \multicolumn{1}{c}{\textbf{ 34.21}}	&  \multicolumn{1}{c}{\textbf{    11.49}} & 
\multicolumn{1}{c}{\textbf{ 27.50}}
\\

\end{tabular}
}

\label{table_lp}
\end{table}

\subsection{Heuristic results}

In order to provide a computational comparison between the work presented in this paper and a standard heuristic from the literature we implemented a time savings based heuristic, where the savings matrix was recalculated as orders were clustered. 
For the estimation of partial route distances
we used~\cite{valle2017}.
As mentioned previously above making use of a time savings based heuristic in order to provide a basis for evaluating the performance of another heuristic is common in the literature, e.g.~see~\cite{albareda2009, deKoster1999, henn2012_2, henn2012, hong2017, matusiak2017, gils18a,  zulj18}. 

Above we presented a partial integer optimisation (PIO) heuristic, a
matheuristic that was directly based upon our mathematical distance approximation formulation. That heuristic successively batches orders for $\tau$ trolleys at a time, until  all orders have been batched.

For space reasons we do not present here detailed results for either the time savings heuristic or our PIO heuristic, rather Table~\ref{table_tau} shows average computation times and average percentage deviations.  For example in that table for the multiple block test problems with $\Delta=5$ the average computation time for our distance approximation  (as taken from Table~\ref{table2}) is 29.7 seconds with the average percentage deviations from the routing results for the optimal approach being 0.01\% for 
no-reversal routing and 1.30\% for reversal routing. Our PIO heuristic with $\tau=1$ has an average computation time of 0.3 seconds with the corresponding average percentage deviations being 1.70\% and 1.17\%. 

The computation time for batching  for our distance approximation  and for our PIO heuristic does not depend on the routing approach (no-reversal or reversal) adopted. However the computation time for the time savings heuristic does depend upon the routing approach adopted.
For the time savings heuristic for multiple block test problems with $\Delta=5$ the average computation time for no-reversal routing was 14.0 seconds, with the average percentage deviation being 13.75\%, whilst for reversal routing the average computation time was 19.3 seconds with the average percentage deviation being 3.92\%. 

Note here that in order to make a consistent comparison all of the routing results shown in Table~\ref{table_tau} were generated by taking the order batching as decided by each of the approaches shown and optimally routing  each individual picker/trolley using~\cite{valle2017}.

Examining 
Table~\ref{table_tau}
we can see that for the single block test problems  our distance approximation   significantly outperforms the time savings heuristic, both with regard to computation time and routing quality (as noted previously above with respect to 
Table~\ref{table_single_extra}). For these instances the average computation times for our distance approximation are very small. 

For the multiple block test problems  our distance approximation   has larger average computation times. For two of the three $\Delta$ cases our distance approximation  has a higher average computation time but significantly better routing quality than the time savings heuristic. However for all three $\Delta$ cases our PIO heuristic (which is itself  directly based on partial integer optimisation of our mathematical distance approximation formulation) significantly outperforms the time savings heuristic, both with regard to computation time and routing quality, for $\tau=1$ and for $\tau=2$. 

For the multiple block test problems,  comparing our PIO heuristic against the distance approximation approach as in Table~\ref{table2} we can see that it has a far smaller computation time, but gives results of reasonable quality. Over all the multiple block $\Delta$ values in Table~\ref{table_tau} the distance approximation approach has an average computation  time of 272.8	seconds with average percentage deviations of $-0.07$\% for no-reversal routing and 	
1.70\% for reversal routing. Our PIO heuristic has an average computation  time of 0.4	seconds with average percentage deviations of 1.63\% for no-reversal routing and 1.58\% for reversal routing for $\tau=1$. For $\tau=2$ the average computational time is  6.1 seconds,	with the percentage deviations being 0.75\% for no-reversal routing and 	2.16\% for reversal routing.

\begin{table}[!ht]
\centering

\caption{Heuristic results}
\vspace{0.4cm}

{\scriptsize
\renewcommand{\tabcolsep}{2.9mm} 
\renewcommand{\arraystretch}{1.1} 
\begin{tabular}{|cccrcc|}
\hline
 &$\Delta$ & Case & T(s) & No-reversal  & Reversal  \\
& & & & routing \% dev. & routing \% dev. \\
\hline
Multiple  & 5 & Distance approximation, Table~\ref{table2} & 29.7& 0.01 & 1.30 \\
blocks & & PIO $\tau=1$ & 0.3 & 1.70 & 1.17     \\
       & & PIO $\tau=2$ & 8.9 & 0.39 & 1.60     \\
       & & PIO $\tau=3$ & 25.6 & 0.14 & 1.35     \\

& & Time savings & 14.0, 10.3 & 13.75 & 3.92 \\

& 10 & Distance approximation, Table~\ref{table2} & 10.7& 0.01 & 1.55 \\
& & PIO $\tau=1$ &   0.3 & 1.59 & 1.16   \\
& & PIO $\tau=2$ &   5.9 & 0.88 & 2.28   \\
& & PIO $\tau=3$ &  22.9 & 0.01 & 1.56   \\

& & Time savings & 23.9, 19.1 & 8.21 & 8.92 \\

& 20 &  Distance approximation, Table~\ref{table2} & 777.8 & -0.23 & 2.26 \\

& & PIO $\tau=1$ &   0.5 & 1.60 & 2.40   \\
& & PIO $\tau=2$ &   3.4 & 0.98 & 2.61   \\
& & PIO $\tau=3$ &  47.0 & 0.12 & 2.83   \\

& & Time savings & 26.7, 23.7& 4.73 & 9.98\\

\hline
Single  & 5 & Distance approximation, Table~\ref{table3} & 0.6& 0 & 2.52 \\
block  & & PIO $\tau=1$ & 0.1 & 1.57 & 3.10     \\
   & & PIO $\tau=2$ & 0.1   & 0.42 & 2.40  \\
       & & PIO $\tau=3$ & 21.55 & 0.11 & 2.42  \\

& & Time savings & 7.1, 7.1 & 3.94 & 5.54 \\

& 10 & Distance approximation, Table~\ref{table3} & 0.1 & 0 & 3.31 \\
& & PIO $\tau=1$ & 0.1 & 0.60 & 1.83     \\
& & PIO $\tau=2$ &   0.1 & 0.40 & 2.60     \\
& & PIO $\tau=3$ & 21.52 & 0.00 & 2.29     \\

& & Time savings & 7.3, 11.9 & 3.87 & 13.57 \\

& 20 & Distance approximation, Table~\ref{table3} & 0.6 & 0 & 2.15\\
& & PIO $\tau=1$ & 0.2 & 0.58 & 3.30     \\
& & PIO $\tau=2$ &   0.2 & 0.31 & 2.45     \\
& & PIO $\tau=3$ & 43.47 & 0.03 & 1.85     \\

& & Time savings & 6.9, 17.3 & 6.38 & 3.15 \\     

\hline
\end{tabular}
}
\label{table_tau}
\end{table}

\subsection{Larger test problems}

The test problems dealt with above had $|O| \leq 30$, with $W_A=8$ and $T \leq 6$. In order to consider some larger test problems we randomly generated some further problems (in the same manner as in~\cite{valle2017}) with $|O|=25,30,50,75$ and 
$W_A=8,16$.

In Table~\ref{table2} and Table~\ref{table3}  we compared our distance approximation approach with the results given 
by the optimal branch-and-cut approach~\cite{valle2017}. 
We also applied that optimal approach  to these larger problems. The solutions obtained are shown in Table~\ref{table7a}.

\begin{table}[!ht]
\centering

\caption{Solution values using an optimal branch-and-cut 
approach
 for the larger problems}
\vspace{0.4cm}

{\scriptsize
\renewcommand{\tabcolsep}{2.9mm} 
\renewcommand{\arraystretch}{1.2} 
\begin{tabular}{|ccc|c|c|c|c|}
\hline
\multirow{3}{*}{$\Delta$} & \multirow{3}{*}{$|O|$} & \multirow{3}{*}{$W_A$} & \multicolumn{2}{c|}{Single block} & \multicolumn{2}{c|}{Multiple blocks}\\
\cline{4-7}
 &  &  & \multicolumn{1}{c|}{No-reversal routing} & \multicolumn{1}{c|}{Reversal routing} & \multicolumn{1}{c|}{No-reversal routing} & \multicolumn{1}{c|}{Reversal routing}\\
\hline
5	&	25	&	8	&	1268.4	&	1180.9	&	1155.1	&	1123.5	\\
 	&	  	&	16	&	1520.2	&	1407.0	&	1514.5	&	1407.9	\\
 	&	30	&	8	&	1358.5	&	1269.1	&	1265.4	&	1263.5	\\
 	&	  	&	16	&	1590.2	&	1493.1	&	1610.9	&	1573.4	\\
 	&	50	&	8	&	2122.3	&	2075.2	&	2134.8	&	2504.1	\\
 	&	  	&	16	&	2544.0	&	2685.2	&	2665.5	&	2906.2	\\
 	&	75	&	8	&	2901.9	&	3473.8	&	2889.0	&	3335.7	\\
 	&	  	&	16	&	3601.6	&	4144.6	&	3747.3	&	11106.3	\\
\hline													
10	&	25	&	8	&	1278.4	&	1200.5	&	1264.5	&	1266.1	\\
 	&	  	&	16	&	1530.2	&	1484.5	&	1635.5	&	1505.6	\\
 	&	30	&	8	&	1350.4	&	1278.5	&	1362.4	&	1345.6	\\
 	&	  	&	16	&	1668.3	&	1573.4	&	1699.4	&	1678.4	\\
 	&	50	&	8	&	2266.3	&	2217.3	&	2211.3	&	2354.5	\\
 	&	  	&	16	&	2886.8	&	3098.1	&	2781.3	&	3209.2	\\
 	&	75	&	8	&	3121.8	&	3653.4	&	3105.2	&	8599.9	\\
 	&	  	&	16	&	3817.4	&	4156.7	&	3915.5	&	11611.4	\\
\hline													
20	&	25	&	8	&	1791.0	&	1697.2	&	1729.2	&	1692.3	\\
 	&	  	&	16	&	2110.9	&	2051.7	&	2187.9	&	2025.1	\\
 	&	30	&	8	&	2069.7	&	1990.6	&	2062.2	&	1944.7	\\
 	&	  	&	16	&	2575.7	&	2390.3	&	2633.0	&	2805.3	\\
 	&	50	&	8	&	2862.8	&	2819.8	&	2801.3	&	7594.5	\\
 	&	  	&	16	&	3338.4	&	9078.4	&	3711.9	&	11169.3	\\
 	&	75	&	8	&	4039.3	&	4577.7	&	3987.2	&	10006.3	\\
 	&	  	&	16	&	4858.7	&	11843.4	&	5086.6	&	15577.7	\\

\hline
\multicolumn{7}{|l|}{Only three of the 96 cases seen above terminated before the 6 hour time limit}  \\
\multicolumn{7}{|l|}{Hence the average computation time for the problems shown  is (effectively) 21600 seconds} \\
\hline
\end{tabular}
}
\label{table7a}
\end{table}

Table~\ref{table_da_sl} shows the results for our distance approximation approach on the single block larger problems. Here the percentage deviation values given are calculated with respect to the values shown in Table~\ref{table7a}. Considering Table~\ref{table_da_sl} we can see from the problems solved to proven optimality that (for a given $\Delta$ and $|O|$) increasing $W_A$ has a very significant effect on computation time. Note however that, as shown by the negative average percentage deviation values, the results obtained are superior to those obtained by applying the approach of~\cite{valle2017}, which required an average of 21600 seconds of computation.

One point to note from Table~\ref{table_da_sl} is that even for those single block instances where we did not solve our distance approximation approach to proven optimality the best-known solution value found at time limit termination corresponded to the subsequent routing solution value.

\begin{table}[!ht]
\centering
\caption{Single block distance approximation results: larger problems
}
\vspace{0.4cm}

{\scriptsize
\renewcommand{\tabcolsep}{0.75mm} 
\renewcommand{\arraystretch}{1.2} 
\begin{tabular}{|cccc|rrrr|rrr|rrr|}
\hline
\multirow{3}{*}{$\Delta$} & \multirow{3}{*}{$|O|$} & \multirow{3}{*}{$W_A$} & \multirow{3}{*}{$T$} & \multicolumn{4}{c|}{Results, distance approximation} & \multicolumn{3}{c|}{No-reversal routing} & \multicolumn{3}{c|}{Reversal routing}\\
 &  &  &  & \multicolumn{1}{c}{T(s)} & \multicolumn{1}{c}{Optimal} & \multicolumn{1}{c}{Root node} & \multicolumn{1}{c|}{Number nodes} & \multicolumn{1}{c}{T(s)} & \multicolumn{1}{c}{Value} & \multicolumn{1}{c|}{\% dev.} & \multicolumn{1}{c}{T(s)} & \multicolumn{1}{c}{Value} & \multicolumn{1}{c|}{\% dev.}\\
 &  &  &  &  & (best-known) & \multicolumn{1}{c}{lower bound} &  &  &  &  &  &  & \\
\hline

5	&	25	&	8	&	4	&	1.7	&	1268.4	&	870.4	&	1922	&	0.3	&	1268.4	&	0	&	1.2	&	1252.9	&	6.10	\\
 	&	  	&	16	&	4	&	72.4	&	1500.2	&	900.1	&	61339	&	1.7	&	1500.2	&	-1.32	&	16.8	&	1457.3	&	3.57	\\
 	&	30	&	8	&	4	&	6.6	&	1358.5	&	1084.4	&	2588	&	0.4	&	1358.5	&	0	&	1.5	&	1297.1	&	2.21	\\
 	&	  	&	16	&	4	&	252.8	&	1590.2	&	990.3	&	87901	&	2.0	&	1590.2	&	0	&	7.6	&	1537.1	&	2.95	\\
 	&	50	&	8	&	7	&	20303.4	&	2070.2	&	1494.7	&	3377840	&	0.3	&	2070.2	&	-2.45	&	1.7	&	2066.5	&	-0.42	\\
 	&	  	&	16	&	7	&	21600.0	&	(2545.6)	&	1365.7	&	1823983	&	7.5	&	2545.6	&	0.06	&	18.4	&	2488.0	&	-7.34	\\
 	&	75	&	8	&	10	&	21600.0	&	(2902.2)	&	2055.9	&	310409	&	0.5	&	2902.2	&	0.01	&	2.2	&	2887.2	&	-16.89	\\
 	&	  	&	16	&	10	&	21600.0	&	(3591.4)	&	1932.5	&	135565	&	6.2	&	3591.4	&	-0.28	&	13.5	&	3500.6	&	-15.54	\\

 &  &  & \textbf{Average:} & \textbf{ 10679.6} &  &  &  & \textbf{     2.4} &  & \textbf{   -0.50} & \textbf{     7.9} &  & \textbf{   -3.17}\\
\hline

10	&	25	&	8	&	4	&	0.6	&	1278.4	&	1038.5	&	535	&	0.1	&	1278.4	&	0	&	2.0	&	1278.4	&	6.49	\\
 	&	  	&	16	&	4	&	17.5	&	1530.2	&	1012.7	&	10092	&	2.7	&	1530.2	&	0	&	21.8	&	1496.4	&	0.80	\\
 	&	30	&	8	&	4	&	1.1	&	1350.4	&	1176.5	&	240	&	0.5	&	1350.4	&	0	&	4.3	&	1318.4	&	3.12	\\
 	&	  	&	16	&	4	&	103.2	&	1628.3	&	1172.9	&	31375	&	2.0	&	1628.3	&	-2.40	&	26.7	&	1602.4	&	1.84	\\
 	&	50	&	8	&	7	&	2081.1	&	2194.3	&	1670.2	&	414421	&	0.4	&	2194.3	&	-3.18	&	2.6	&	2158.3	&	-2.66	\\
 	&	  	&	16	&	7	&	21600.0	&	(2702.1)	&	1570.6	&	1825500	&	38.9	&	2702.1	&	-6.40	&	16.7	&	2630.8	&	-15.08	\\
 	&	75	&	8	&	10	&	21600.0	&	(3111.8)	&	2314.3	&	293148	&	0.6	&	3111.8	&	-0.32	&	3.2	&	3000.4	&	-17.87	\\
 	&	  	&	16	&	10	&	21600.0	&	(3881.6)	&	2083.8	&	174119	&	12.9	&	3881.6	&	1.68	&	39.5	&	3754.9	&	-9.67	\\

 &  &  & \textbf{Average:} & \textbf{  8375.4} &  &  &  & \textbf{     7.3} &  & \textbf{   -1.33} & \textbf{    14.6} &  & \textbf{   -4.13}\\
\hline

20	&	25	&	8	&	5	&	0.3	&	1791.0	&	1647.0	&	198	&	0.4	&	1791.0	&	0	&	2.7	&	1729.6	&	1.91	\\	
 	&	  	&	16	&	5	&	103.2	&	2062.7	&	1668.8	&	50646	&	4.0	&	2062.7	&	-2.28	&	22.9	&	2028.8	&	-1.12	\\	
 	&	30	&	8	&	6	&	8.0	&	2069.7	&	1759.5	&	4504	&	0.3	&	2069.7	&	0	&	3.0	&	2016.2	&	1.29	\\	
 	&	  	&	16	&	6	&	2854.4	&	2393.3	&	1786.6	&	923747	&	3.4	&	2393.3	&	-7.08	&	15.4	&	2347.2	&	-1.80	\\	
 	&	50	&	8	&	8	&	3376.0	&	2862.8	&	2538.3	&	1646129	&	0.7	&	2862.8	&	0	&	8.0	&	2763.4	&	-2.00	\\	
 	&	  	&	16	&	8	&	21600.0	&	(3336.2)	&	2353.8	&	639722	&	4.8	&	3336.2	&	-0.07	&	27.3	&	3250.0	&	-64.20	\\	
 	&	75	&	8	&	12	&	21600.0	&	(4039.3)	&	3336.8	&	223514	&	1.3	&	4039.3	&	0	&	7.9	&	3894.1	&	-14.93	\\	
 	&	  	&	16	&	12	&	21600.0	&	(4706.8)	&	3028.3	&	68850	&	10.5	&	4706.8	&	-3.13	&	29.7	&	4570.0	&	-61.41	\\

 &  &  & \textbf{Average:} & \textbf{  8892.7} &  &  &  & \textbf{     3.2} &  & \textbf{   -1.57} & \textbf{    14.6} &  & \textbf{  -17.78}\\
\hline
\multicolumn{4}{r}{\textbf{Average all:}} & \textbf{  9315.9} &  &  & \multicolumn{1}{c}{} & \textbf{     4.3} &  & \multicolumn{1}{c}{\textbf{   -1.13}} & \textbf{    12.4} &  & \multicolumn{1}{c}{\textbf{   -8.36}}
\end{tabular}
}
\label{table_da_sl}
\end{table}

It is clear from Table~\ref{table_da_sl} that as problem size increases finding the exact (optimal) solution to our distance approximation formulation becomes increasingly difficult. For this reason we might, for larger problems, utilise our partial integer optimisation (PIO) heuristic, which is itself directly based upon our mathematical distance approximation formulation.

Table~\ref{large_h_s} shows the results for our PIO heuristic and the time savings heuristic on the single block test problems considered in Table~\ref{table_da_sl}. Here our PIO heuristic (on average) performs significantly better than the time savings heuristic, both in terms of computation time and in terms of quality of solution. As compared with Table~\ref{table_da_sl} our PIO heuristic does not (on average) produce results of the same quality, but is significantly faster.

Over all the single block $\Delta$ values in Table~\ref{table_da_sl} and Table~\ref{large_h_s}
the distance approximation approach has an average computation  time of 9315.9 	seconds  (ignoring the times required for routing the order batching decided for convenience of comparison) with average percentage deviations of $-1.13$\% for no-reversal routing 
and $-8.36$\% for
reversal routing. Our PIO heuristic has an average computation  time of 196.5	seconds with average percentage deviations of 2.41\% for no-reversal routing and 
$-6.36$\% for reversal routing.

\begin{table}[!ht]
\centering

\caption{Single block heuristic results: larger problems
}
\vspace{0.4cm}

{\scriptsize
\renewcommand{\tabcolsep}{1.75mm} 
\renewcommand{\arraystretch}{1.2} 
\begin{tabular}{|ccc|r|rr|rr|rrr|rrr|}
\hline
\multirow{3}{*}{$\Delta$} & \multirow{3}{*}{$|O|$} & \multirow{3}{*}{$W_A$} & \multicolumn{5}{c|}{PIO $\tau = 1$} & \multicolumn{6}{c|}{Time savings}\\
 &  &  & \multicolumn{1}{c}{} & \multicolumn{2}{c}{No-reversal} & \multicolumn{2}{c|}{Reversal} & \multicolumn{3}{c}{No-reversal} & \multicolumn{3}{c|}{Reversal}\\
 &  &  & \multicolumn{1}{c}{T(s)} & \multicolumn{1}{c}{Value} & \multicolumn{1}{c}{\% dev.} & \multicolumn{1}{c}{Value} & \multicolumn{1}{c|}{\% dev.} & \multicolumn{1}{c}{T(s)} & \multicolumn{1}{c}{Value} & \multicolumn{1}{c}{\% dev.} & \multicolumn{1}{c}{T(s)} & \multicolumn{1}{c}{Value} & \multicolumn{1}{c|}{\% dev.}\\
\hline

5	&	25	&	8	&	0.3	&	1350.4	&	6.46	&	1283.0	&	8.65	&	18.8	&	1350.4	&	6.46	&	17.1	&	1278.6	&	8.27	\\
 	&	  	&	16	&	1.0	&	1520.2	&	0	&	1424.8	&	1.27	&	229.8	&	1720.2	&	13.16	&	113.2	&	1478.7	&	5.10	\\
 	&	30	&	8	&	1.0	&	1432.4	&	5.44	&	1378.9	&	8.65	&	27.0	&	1627.1	&	19.77	&	23.5	&	1508.3	&	18.85	\\
 	&	  	&	16	&	1.8	&	1640.2	&	3.14	&	1503.4	&	0.69	&	319.2	&	1770.2	&	11.32	&	157.0	&	1563.2	&	4.69	\\
 	&	50	&	8	&	10.9	&	2202.4	&	3.77	&	2112.4	&	1.79	&	74.1	&	2348.3	&	10.65	&	57.3	&	2190.3	&	5.55	\\
 	&	  	&	16	&	18.9	&	2683.7	&	5.49	&	2542.4	&	-5.32	&	739.3	&	2913.7	&	14.53	&	425.1	&	2583.5	&	-3.79	\\
 	&	75	&	8	&	410.7	&	3114.0	&	7.31	&	2936.1	&	-15.48	&	159.4	&	3326.2	&	14.62	&	110.2	&	3130.2	&	-9.89	\\
 	&	  	&	16	&	179.9	&	3713.2	&	3.10	&	3576.6	&	-13.70	&	1456.6	&	3981.2	&	10.54	&	891.5	&	3614.8	&	-12.78	\\

 &  & \textbf{Average:} & \textbf{    78.1} &  & \textbf{     4.34} &  & \textbf{    -1.68} & \textbf{   378.0} &  & \textbf{    12.63} & \textbf{   224.4} &  & \textbf{     2.0}\\
\hline
10	&	25	&	8	&	0.5	&	1278.4	&	0	&	1276.6	&	6.34	&	19.5	&	1350.4	&	5.63	&	25.9	&	1355.1	&	12.88	\\
 	&	  	&	16	&	1.2	&	1560.2	&	1.96	&	1496.9	&	0.84	&	323.9	&	1730.2	&	13.07	&	162.8	&	1594.9	&	7.44	\\
 	&	30	&	8	&	0.6	&	1350.4	&	0	&	1314.4	&	2.81	&	28.7	&	1360.4	&	0.74	&	34.5	&	1625.4	&	27.13	\\
 	&	  	&	16	&	3.1	&	1670.2	&	0.11	&	1652.2	&	5.01	&	422.7	&	2010.8	&	20.53	&	204.7	&	1901.5	&	20.85	\\
 	&	50	&	8	&	8.5	&	2360.2	&	4.14	&	2292.6	&	3.40	&	79.8	&	2440.2	&	7.67	&	79.2	&	2225.4	&	0.37	\\
 	&	  	&	16	&	26.0	&	2845.9	&	-1.42	&	2690.8	&	-13.15	&	984.2	&	3005.6	&	4.12	&	495.7	&	2677.9	&	-13.56	\\
 	&	75	&	8	&	675.9	&	3215.9	&	3.01	&	3146.9	&	-13.86	&	176.2	&	3348.0	&	7.25	&	153.5	&	3220.5	&	-11.85	\\
 	&	  	&	16	&	188.8	&	3849.4	&	0.84	&	3726.2	&	-10.36	&	1980.6	&	4281.0	&	12.14	&	1037.6	&	3845.0	&	-7.50	\\

 &  & \textbf{Average:} & \textbf{   113.1} &  & \textbf{     1.08} &  & \textbf{    -2.37} & \textbf{   502.0} &  & \textbf{     8.89} & \textbf{   274.2} &  & \textbf{     4.47}\\
\hline
20	&	25	&	8	&	0.5	&	1803.0	&	0.67	&	1799.0	&	6.00	&	19.2	&	1801.0	&	0.56	&	41.7	&	1829.0	&	7.77	\\
 	&	  	&	16	&	1.3	&	2192.7	&	3.88	&	2093.6	&	2.04	&	366.2	&	2202.7	&	4.35	&	195.5	&	2020.2	&	-1.54	\\
 	&	30	&	8	&	1.1	&	2151.6	&	3.96	&	2065.7	&	3.77	&	27.4	&	2153.6	&	4.05	&	57.1	&	2041.6	&	2.56	\\
 	&	  	&	16	&	4.8	&	2525.2	&	-1.96	&	2428.4	&	1.59	&	536.8	&	2583.3	&	0.30	&	273.4	&	2370.8	&	-0.82	\\

 	&	50	&	8	&	1487.0	&	2874.7	&	0.42	&	2838.8	&	0.67	&	76.2	&	3223.4	&	12.60	&	136.9	&	2884.8	&	2.31	\\
 	&	  	&	16	&	45.1	&	3490.3	&	4.55	&	3391.2	&	-62.65	&	1504.2	&	3891.0	&	16.55	&	747.1	&	3607.9	&	-60.26	\\
 	&	75	&	8	&	1461.1	&	4161.2	&	3.02	&	4057.3	&	-11.37	&	173.5	&	4297.1	&	6.38	&	266.3	&	4085.6	&	-10.75	\\
 	&	  	&	16	&	186.1	&	4856.7	&	-0.04	&	4711.8	&	-60.22	&	3054.5	&	5098.6	&	4.94	&	1464.2	&	4748.4	&	-59.91	\\

 &  & \textbf{Average:} & \textbf{   398.4} &  & \textbf{     1.81} &  & \textbf{   -15.02} & \textbf{   719.8} &  & \textbf{     6.22} & \textbf{   397.8} &  & \textbf{   -15.08}\\
\hline
\multicolumn{3}{r}{\textbf{Average all:}} & \multicolumn{1}{c}{\textbf{   196.5}} &  & \multicolumn{1}{c}{\textbf{     2.41}} &  & \multicolumn{1}{c}{\textbf{    -6.36}} & \multicolumn{1}{c}{\textbf{   533.2}} &  & \multicolumn{1}{c}{\textbf{     9.25}} & \multicolumn{1}{c}{\textbf{   298.8}} &  & \multicolumn{1}{c}{\textbf{    -2.87}}\\
\end{tabular}
}

\label{large_h_s}
\end{table}

Table~\ref{table_da_ml} shows the results for our distance approximation approach on the multiple block larger problems. As for Table~\ref{table_da_sl} we can see from the problems solved to proven optimality that (for a given $\Delta$ and $|O|$) increasing $W_A$ has a very significant effect on computation time. Also, as for Table~\ref{table_da_sl},  the results obtained are superior to those obtained by applying the approach 
of~\cite{valle2017}. In particular note here the very high negative percentage deviations that are seen for reversal routing, indicating that our distance approximation approach is providing solutions of much higher quality than the approach of~\cite{valle2017}.

Above in relation to the multiple block instances in 
Table~\ref{table2} we found that the values of the optimal distance approximation solution and the routing solution based on the order batching given by the distance approximation were effectively identical (an average deviation of 0.03\% for no-reversal routing). For the larger multiple block instances shown in 
Table~\ref{table_da_ml} the average value of 100(Routing solution value - Optimal or best-known, distance approximation value)/(Optimal, or best-known, distance approximation value) is 0.06\% for no-reversal routing. Hence we can conclude here that for multiple block no-reversal routing (for the test problems examined) our distance approximation very closely approximates the underlying routing distance.

\begin{table}[!ht]
\centering

\caption{Multiple block distance approximation results: larger problems
}
\vspace{0.4cm}

{\scriptsize
\renewcommand{\tabcolsep}{0.75mm} 
\renewcommand{\arraystretch}{1.2} 
\begin{tabular}{|cccc|rrrr|rrr|rrr|}
\hline
\multirow{3}{*}{$\Delta$} & \multirow{3}{*}{$|O|$} & \multirow{3}{*}{$W_A$} & \multirow{3}{*}{$T$} & \multicolumn{4}{c|}{Results, distance approximation} & \multicolumn{3}{c|}{No-reversal routing} & \multicolumn{3}{c|}{Reversal routing}\\
 &  &  &  & \multicolumn{1}{c}{T(s)} & \multicolumn{1}{c}{Optimal} & \multicolumn{1}{c}{Root node} & \multicolumn{1}{c|}{Number nodes} & \multicolumn{1}{c}{T(s)} & \multicolumn{1}{c}{Value} & \multicolumn{1}{c|}{\% dev.} & \multicolumn{1}{c}{T(s)} & \multicolumn{1}{c}{Value} & \multicolumn{1}{c|}{\% dev.}\\
 &  &  &  &  & (best-known)  & \multicolumn{1}{c}{lower bound} &  &  &  &  &  &  & \\
\hline

5	&	25	&	8	&	4	&	30.1	&	1154.2	&	721.7	&	22221	&	0.2	&	1155.1	&	0	&	2.1	&	1131.6	&	0.72	\\
 	&	  	&	16	&	4	&	1616.9	&	1444.0	&	854.1	&	649360	&	4.4	&	1445.3	&	-4.57	&	47.3	&	1355.3	&	-3.74	\\
 	&	30	&	8	&	4	&	384.5	&	1264.4	&	802.1	&	150557	&	0.2	&	1265.4	&	0	&	1.9	&	1219.9	&	-3.45	\\
 	&	  	&	16	&	4	&	21600.0	&	(1578.3)	&	896.8	&	3207670	&	3.4	&	1579.7	&	-1.93	&	37.5	&	1489.4	&	-5.34	\\
 	&	50	&	8	&	7	&	21600.0	&	(2124.3)	&	1070.2	&	1780350	&	0.8	&	2125.2	&	-0.45	&	2.4	&	1998.0	&	-20.21	\\
 	&	  	&	16	&	7	&	21600.0	&	(2585.8)	&	1023.1	&	1369222	&	19.1	&	2588.4	&	-2.89	&	47.4	&	2404.1	&	-17.28	\\
 	&	75	&	8	&	10	&	21600.0	&	(2952.1)	&	1466.1	&	144040	&	0.8	&	2954.1	&	2.25	&	3.9	&	2848.2	&	-14.61	\\
 	&	  	&	16	&	10	&	21600.0	&	(3693.6)	&	1324.6	&	254133	&	65.0	&	3696.3	&	-1.36	&	65.7	&	3413.4	&	-69.27	\\

 &  &  & \textbf{Average:} & \textbf{ 13753.9} &  &  &  & \textbf{    11.7} &  & \textbf{   -1.12} & \textbf{    26.0} &  & \textbf{  -16.65}\\
\hline

10	&	25	&	8	&	4	&	11.5	&	1264.4	&	868.2	&	7276	&	0.9	&	1264.6	&	0.01	&	2.7	&	1233.7	&	-2.56	\\
 	&	  	&	16	&	4	&	4304.8	&	1562.3	&	927.5	&	2245139	&	7.8	&	1564.1	&	-4.37	&	200.6	&	1478.0	&	-1.83	\\
 	&	30	&	8	&	4	&	137.5	&	1362.4	&	975.5	&	39598	&	0.4	&	1362.9	&	0.03	&	4.9	&	1314.1	&	-2.34	\\
 	&	  	&	16	&	4	&	21600.0	&	(1654.2)	&	978.3	&	4204508	&	7.0	&	1655.6	&	-2.58	&	98.1	&	1548.1	&	-7.76	\\
 	&	50	&	8	&	7	&	21600.0	&	(2180.3)	&	1277.8	&	2095008	&	0.4	&	2181.4	&	-1.35	&	7.5	&	2121.3	&	-9.91	\\
 	&	  	&	16	&	7	&	21600.0	&	(2720.0)	&	1200.3	&	1185330	&	30.3	&	2722.2	&	-2.12	&	131.8	&	2584.3	&	-19.47	\\
 	&	75	&	8	&	10	&	21600.0	&	(3093.8)	&	1746.5	&	146444	&	1.0	&	3096.7	&	-0.27	&	6.9	&	3022.2	&	-64.86	\\
 	&	  	&	16	&	10	&	21600.0	&	(3861.6)	&	1581.2	&	149009	&	15.9	&	3863.5	&	-1.33	&	79.9	&	3593.6	&	-69.05	\\

 &  &  & \textbf{Average:} & \textbf{ 14056.7} &  &  &  & \textbf{     8.0} &  & \textbf{   -1.50} & \textbf{    66.6} &  & \textbf{  -22.22}\\
\hline

20	&	25	&	8	&	5	&	221.8	&	1727.0	&	1289.6	&	118215	&	0.3	&	1727.2	&	-0.11	&	3.4	&	1672.4	&	-1.17	\\
 	&	  	&	16	&	5	&	21600.0	&	(2098.7)	&	1311.0	&	3584178	&	12.4	&	2100.1	&	-4.01	&	73.6	&	1966.3	&	-2.90	\\
 	&	30	&	8	&	6	&	10639.7	&	1995.4	&	1344.0	&	3808495	&	0.5	&	1996.2	&	-3.20	&	4.8	&	1919.1	&	-1.32	\\
 	&	  	&	16	&	6	&	21600.0	&	(2445.3)	&	1424.6	&	3220834	&	50.1	&	2447.1	&	-7.06	&	136.3	&	2275.6	&	-18.88	\\
 	&	50	&	8	&	8	&	21600.0	&	(2762.7)	&	1906.8	&	638213	&	0.7	&	2764.2	&	-1.32	&	8.3	&	2729.2	&	-64.06	\\
 	&	  	&	16	&	8	&	21600.0	&	(3438.3)	&	1818.0	&	409990	&	19.5	&	3440.6	&	-7.31	&	115.2	&	3230.0	&	-71.08	\\
 	&	75	&	8	&	12	&	21600.0	&	(3897.3)	&	2444.1	&	54912	&	1.4	&	3899.4	&	-2.20	&	10.5	&	3786.2	&	-62.16	\\
 	&	  	&	16	&	12	&	21600.0	&	(4762.6)	&	2283.9	&	60168	&	46.1	&	4764.5	&	-6.33	&	167.5	&	4517.4	&	-71.00	\\

 &  &  & \textbf{Average:} & \textbf{ 17557.7} &  &  &  & \textbf{    16.4} &  & \textbf{   -3.94} & \textbf{    65.0} &  & \textbf{  -36.57}\\
\hline
\multicolumn{4}{r}{\textbf{Average all:}} & \textbf{ 15122.8} &  &  & \multicolumn{1}{c}{} & \textbf{    12.0} &  & \multicolumn{1}{c}{\textbf{   -2.19}} & \textbf{    52.5} &  & \multicolumn{1}{c}{\textbf{  -25.15}}\end{tabular}
}
\label{table_da_ml}
\end{table}

Table~\ref{large_h_m} shows the results for our partial integer optimisation (PIO) heuristic and the time savings heuristic on the multiple block test problems considered in Table~\ref{table_da_ml}. Here our PIO heuristic (on average) again performs significantly better than the time savings heuristic, both in terms of computation time and in terms of quality of solution. As compared with Table~\ref{table_da_ml} our PIO heuristic does not (on average) produce results of the same quality, but is significantly faster. 

Over all the multiple block $\Delta$ values in Table~\ref{table_da_ml} and Table~\ref{large_h_m}
the distance approximation approach has an average computation  time of 15122.8	seconds  (ignoring the times required for routing the order batching decided for convenience of comparison) with average percentage deviations of $-2.19$\% for no-reversal routing and 
$-25.15$\% for reversal routing. Our PIO heuristic has an average computation  time of 425.6	seconds with average percentage deviations of 0.65\% for no-reversal routing and $-23.51$\% for reversal routing.

\begin{table}[!ht]
\centering

\caption{Multiple block heuristic results: larger problems
}
\vspace{0.4cm}

{\scriptsize
\renewcommand{\tabcolsep}{1.75mm} 
\renewcommand{\arraystretch}{1.2} 

\begin{tabular}{|ccc|r|rr|rr|rrr|rrr|}
\hline
\multirow{3}{*}{$\Delta$} & \multirow{3}{*}{$|O|$} & \multirow{3}{*}{$W_A$} & \multicolumn{5}{c|}{PIO $\tau = 1$} & \multicolumn{6}{c|}{Time savings}\\
 &  &  & \multicolumn{1}{c}{} & \multicolumn{2}{c}{No-reversal} & \multicolumn{2}{c|}{Reversal} & \multicolumn{3}{c}{No-reversal} & \multicolumn{3}{c|}{Reversal}\\
 &  &  & \multicolumn{1}{c}{T(s)} & \multicolumn{1}{c}{Value} & \multicolumn{1}{c}{\% dev.} & \multicolumn{1}{c}{Value} & \multicolumn{1}{c|}{\% dev.} & \multicolumn{1}{c}{T(s)} & \multicolumn{1}{c}{Value} & \multicolumn{1}{c}{\% dev.} & \multicolumn{1}{c}{T(s)} & \multicolumn{1}{c}{Value} & \multicolumn{1}{c|}{\% dev.}\\
\hline
5	&	25	&	8	&	1.2	&	1226.9	&	6.22	&	1137.8	&	1.27	&	39.2	&	1450.4	&	25.56	&	24.8	&	1165.9	&	3.77	\\
 	&	  	&	16	&	5.4	&	1514.4	&	-0.01	&	1411.6	&	0.26	&	652.9	&	1614.7	&	6.62	&	400.1	&	1409.7	&	0.13	\\
 	&	30	&	8	&	2.8	&	1332.5	&	5.30	&	1267.9	&	0.35	&	54.6	&	1607.2	&	27.01	&	33.7	&	1414.6	&	11.96	\\
 	&	  	&	16	&	20.5	&	1655.1	&	2.74	&	1516.1	&	-3.64	&	908.1	&	1957.0	&	21.48	&	548.5	&	1508.4	&	-4.13	\\
 	&	50	&	8	&	22.7	&	2098.0	&	-1.72	&	2010.5	&	-19.71	&	145.2	&	2474.4	&	15.91	&	83.1	&	2062.9	&	-17.62	\\
 	&	  	&	16	&	492.5	&	2645.5	&	-0.75	&	2442.7	&	-15.95	&	2302.2	&	2845.0	&	6.73	&	1471.2	&	2484.9	&	-14.50	\\
 	&	75	&	8	&	647.8	&	2993.2	&	3.61	&	2884.8	&	-13.52	&	304.2	&	3473.2	&	20.22	&	163.8	&	2923.0	&	-12.37	\\
 	&	  	&	16	&	1630.5	&	3665.8	&	-2.17	&	3372.0	&	-69.64	&	5036.6	&	3934.7	&	5.00	&	3093.5	&	3402.1	&	-69.37	\\

 &  & \textbf{Average:} & \textbf{   352.9} &  & \textbf{     1.65} &  & \textbf{   -15.07} & \textbf{  1180.4} &  & \textbf{    16.07} & \textbf{   727.3} &  & \textbf{   -12.77}\\
\hline

10	&	25	&	8	&	0.7	&	1302.5	&	3.01	&	1251.3	&	-1.17	&	56.3	&	1488.4	&	17.71	&	40.1	&	1325.1	&	4.66	\\
 	&	  	&	16	&	3.5	&	1607.0	&	-1.74	&	1501.1	&	-0.30	&	790.0	&	1723.3	&	5.37	&	648.6	&	1539.5	&	2.25	\\
 	&	30	&	8	&	3.4	&	1412.5	&	3.68	&	1364.9	&	1.43	&	77.3	&	1687.0	&	23.83	&	49.8	&	1574.5	&	17.01	\\
 	&	  	&	16	&	13.9	&	1689.5	&	-0.58	&	1604.1	&	-4.43	&	1101.9	&	2070.2	&	21.82	&	752.1	&	1860.9	&	10.87	\\
 	&	50	&	8	&	29.0	&	2291.2	&	3.61	&	2184.4	&	-7.22	&	181.4	&	2564.7	&	15.98	&	115.9	&	2180.4	&	-7.39	\\
 	&	  	&	16	&	312.7	&	2771.8	&	-0.34	&	2564.9	&	-20.08	&	2791.7	&	3050.2	&	9.67	&	1682.9	&	2652.3	&	-17.35	\\
 	&	75	&	8	&	1903.4	&	3204.9	&	3.21	&	3125.8	&	-63.65	&	368.8	&	3638.4	&	17.17	&	227.9	&	3131.8	&	-63.58	\\
 	&	  	&	16	&	1540.7	&	3930.2	&	0.38	&	3695.0	&	-68.18	&	5926.8	&	4139.4	&	5.72	&	3530.0	&	3720.9	&	-67.95	\\

 &  & \textbf{Average:} & \textbf{   475.9} &  & \textbf{     1.40} &  & \textbf{   -20.45} & \textbf{  1411.8} &  & \textbf{    14.66} & \textbf{   880.9} &  & \textbf{   -15.19}\\
\hline
20	&	25	&	8	&	1.6	&	1807.2	&	4.51	&	1708.8	&	0.98	&	69.4	&	1862.8	&	7.73	&	53.9	&	1762.9	&	4.17	\\
 	&	  	&	16	&	11.7	&	2189.8	&	0.09	&	2052.5	&	1.35	&	864.1	&	2197.7	&	0.45	&	618.2	&	1988.2	&	-1.82	\\
 	&	30	&	8	&	3.0	&	2113.8	&	2.50	&	1954.8	&	0.52	&	96.5	&	2307.5	&	11.90	&	70.9	&	2021.5	&	3.95	\\
 	&	  	&	16	&	40.4	&	2508.5	&	-4.73	&	2337.9	&	-16.66	&	1329.2	&	2632.2	&	-0.03	&	792.4	&	2311.4	&	-17.61	\\
 	&	50	&	8	&	54.0	&	2851.7	&	1.80	&	2754.6	&	-63.73	&	257.0	&	3347.9	&	19.51	&	169.0	&	3000.0	&	-60.50	\\
 	&	  	&	16	&	714.7	&	3445.3	&	-7.18	&	3277.3	&	-70.66	&	3547.5	&	4042.2	&	8.90	&	2194.6	&	3474.6	&	-68.89	\\
 	&	75	&	8	&	659.4	&	3956.7	&	-0.76	&	3849.0	&	-61.53	&	511.7	&	4463.3	&	11.94	&	326.1	&	3963.4	&	-60.39	\\
 	&	  	&	16	&	2099.1	&	4825.1	&	-5.14	&	4618.9	&	-70.35	&	7145.4	&	5218.6	&	2.60	&	4433.9	&	4649.7	&	-70.15	\\

 &  & \textbf{Average:} & \textbf{   448.0} &  & \textbf{    -1.11} &  & \textbf{   -35.01} & \textbf{  1727.6} &  & \textbf{     7.88} & \textbf{  1082.4} &  & \textbf{   -33.91}\\
\hline
\multicolumn{3}{r}{\textbf{Average all:}} & \multicolumn{1}{c}{\textbf{   425.6}} &  & \multicolumn{1}{c}{\textbf{     0.65}} &  & \multicolumn{1}{c}{\textbf{   -23.51}} & \multicolumn{1}{c}{\textbf{  1439.9}} &  & \multicolumn{1}{c}{\textbf{    12.87}} & \multicolumn{1}{c}{\textbf{   896.9}} &  & \multicolumn{1}{c}{\textbf{   -20.62}}\\
\end{tabular}
}

\label{large_h_m}
\end{table}

\subsection{Computational conclusions}
Given the computational results presented above it seems reasonable to conclude that:
\begin{itemize}
\item the distance approximation approach given in this paper gives results nearly as good, or better, than the approach given in our previous paper~\cite{valle2017}, but in significantly lower computation times (Tables~\ref{table2},\ref{table3},\ref{table7a},\ref{table_da_sl},\ref{table_da_ml})
\item for no-reversal routing our distance
approximation very closely approximates the underlying routing distance (Tables~\ref{table2},\ref{table3},\ref{table_da_sl},\ref{table_da_ml})
\item the symmetry breaking constraints presented are of computational benefit (Tables~\ref{table4},\ref{table5})
\item the constraints  (valid inequalities) that we have presented significantly  improve the value of the linear programming relaxation of our formulation (Table~\ref{table_lp})
\item the partial integer optimisation (PIO) heuristic presented, which is directly based on our mathematical distance approximation formulation, produces good quality results in reasonable computation times, and is superior to a time savings heuristic 
(Tables~\ref{table_tau},\ref{large_h_s},\ref{large_h_m})
\end{itemize}

\section{Conclusions}
\label{sec:conclusions}

In this paper we dealt with the problem of order batching for picker routing. We presented an approach that directly addressed order batching, but which used a distance approximation to influence the batching of orders without directly addressing the routing problem.

We presented a basic formulation based on deciding the orders to be batched together so as to optimise an objective that approximates the picker routing distance travelled (both the distance travelled within aisles and the distance travelled within cross-aisles).  We then extended our formulation 
for cases where we have more than one block in the warehouse. We presented constraints to remove symmetry in order to lessen the computational effort required, as well as  constraints that significantly  improve the value of the linear programming relaxation of our formulation. 
We proved that in the case of a single block with 
no-reversal routing our distance approximation gives an optimal
allocation of orders to pickers and an optimal routing.
A heuristic algorithm based on partial integer optimisation of our formulation was also presented.  Once order batching had been decided we optimally routed each individual picker.

Extensive computational results for publicly available test problems involving up to 75 orders were given for both single and multiple block warehouse configurations. 

In terms of future research we plan to develop a variety of new algorithms for order batching based on  metaheuristic approaches that make use of the approximation presented in this paper to estimate the underlying routing distance. Being able to estimate this distance, without the computational expense of computing picker routes, offers the potential for such algorithms to be better able to deal with large problems. It is  currently an open question as to whether using a distance approximation, rather than an explicit routing, will lead to better overall results for metaheuristic approaches to order batching.

 \clearpage
\newpage
 \pagestyle{empty}
\linespread{1}
\small \normalsize 


\section*{\textbf{Acknowledgments}}
The authors would like to thank the anonymous reviewers for their comments and suggestions to improve the quality of the paper.

\bibliographystyle{plainnat}

\bibliography{tesco}

\end{document}